\documentclass[final]{svjour3}
\usepackage{graphicx}
\usepackage{style}
\usepackage{bm}
\usepackage{epstopdf, epsfig}

\title{An accurate integral equation method for Stokes flow with piecewise smooth boundaries}

\author{Lukas Bystricky \and Sara P\aa lsson \and Anna-Karin Tornberg}
\institute{Department of Mathematics, KTH Royal Institute of Technology, Stockholm, SE\\ \email{lukasby@kth.se}}

\begin{document}

\maketitle

\begin{abstract}

Two-dimensional Stokes flow through a periodic channel is considered. The channel walls need only be Lipschitz continuous, in other words they are allowed to have corners. Boundary integral methods are an attractive tool for numerically solving the Stokes equations, as the partial differential equation can be reformulated into an integral equation that must be solved only over the boundary of the domain. When the boundary is at least $C^1$ smooth, the boundary integral kernel is a compact operator, and traditional Nystr\"{o}m methods can be used to obtain highly accurate solutions. In the case of Lipschitz continuous boundaries, however, obtaining accurate solutions using the standard Nystr\"{o}m method can require high resolution. We adapt a technique known as recursively compressed inverse preconditioning to accurately solve the Stokes equations without requiring any more resolution than is needed to resolve the boundary. Combined with a periodic fast summation method we construct a method that is $\mathcal{O}(N\log N)$ where $N$ is the number of quadrature points on the boundary. We demonstrate the robustness of this method by extending an existing boundary integral method for viscous drops to handle the movement of drops near corners. 

\end{abstract}

\section{Introduction}

The Stokes equations are used to model slowly moving, highly viscous, fluid flow. They can be thought of as the zero Reynolds number limit of the Navier-Stokes equations. Of the many applications of the Stokes equations, they are often used to model particle suspensions including solid particles \cite{Bystricky2019,Klinteberg2016}, drops \cite{Ojala2015,Palsson2019a,Palsson2019b}, or vesicles \cite{Quaife2014}.  Even for non-zero Reynolds numbers, the Stokes equations often describe very well the flow near a solid boundary, and can therefore be used to derive effective slip boundary conditions for problems at higher Reynolds numbers \cite{Achdou1998,Dalibard2011}.

An advantage of using the Stokes equations over the Navier-Stokes equations is that the Stokes equations are linear and elliptic, allowing us to recast them as a boundary integral equation (BIE) \cite{Kim2005,Ladyzhenskaya1987,Pozrikidis1992}. BIEs have several nice properties. All the information needed to solve a BIE is confined to the boundary of the domain; this leads to an immediate dimension reduction. If the boundary of the domain is at least $C^1$ smooth, the Stokes equations can be represented as a second-kind Fredholm equation \cite{Power1987}. After discretization the condition number of the resulting linear system is independent of the number of discretization points used, meaning that very highly accurate solutions are obtainable. Traditionally one major drawback of BIEs was the need to solve dense linear systems. However, by using efficient iterative solvers such as GMRES \cite{Saad1986}, combined with fast matrix vector products \cite{Greengard1987,Lindbo2011}, the cost to solve the $N\times N$ dense linear system can be reduced to $\mathcal{O}(N)$ or $\mathcal{O}(N\log N)$, where $N$ is the number of discretization points on the boundary of the domain.

In this paper, we will be considering wall-bounded, periodic Stokes flow. For particulate flows, such models are useful because they allow for the computation of various time averaged quantities over a relatively small reference cell, without the need to simulate an unfeasibly large domain. BIEs have been successfully applied to such problems \cite{Marple2015,Palsson2019b,Zhao2010}.

Lipshitz continuous boundaries admit corners, specifically they allow for countable number of corners with the angle $\theta$ of each corner satisfying $0 < \theta < 2\pi$ (thus prohibiting cusps). Solving PDEs on Lipschitz domains to high accuracy everywhere in the domain is in general quite challenging, independent of the numerical method used; see for example the discussion in  \cite{Gopal2019}. When solving boundary integral equations on domains with corners, the standard Nystr\"{o}m method fails to achieve optimal accuracy \cite{Bremer2012}. While the underlying equation has a unique solution, the layer density defined on the boundary can become weakly singular at corner points, thus reducing the accuracy of regular quadrature rules, such as composite Gauss-Legendre quadrature.  One approach to solving this problem is to cluster additional quadrature points near the corners. This of course dramatically increases the number of unknowns, while at the same time the accuracy of such an approach is limited \cite{Bremer2012,Helsing2013b}. A recent paper \cite{Wu2019} demonstrates an approach that automates both the spatial adaptivity and the order of the quadrature (similar to $hp$ adaptivity) to achieve high accuracy for complicated domains containing several corners. Other approaches have been developed, some of which use elegant kernel-dependent custom quadrature rules \cite{Rachh2017,Serkh2016}. In this paper, we use a kernel-independent method known as \emph{recursively compressed inverse preconditioning} \cite{Helsing2018,Helsing2013b,Helsing2008a}. This method has the advantage of being relatively simple to implement on top of existing code, and does not introduce any additional unknowns to the linear system that must be solved. We demonstrate the robustness of this method when applied to periodic Stokes flow by solving problems involving moving viscous drops.

\section{Governing Equations}\label{sec:gov}

The governing equations are the steady incompressible Stokes equations,
\begin{subequations}\label{eq:stokes}
\begin{alignat}{2}
	-\mu\Delta \uu(\xx) + \nabla p(\xx) &= \mathbf{0},\qquad & \xx\in \Omega,\\
	\nabla\cdot\uu(\xx) &= 0, \qquad& \xx\in\Omega,
\end{alignat}
\end{subequations}
where $\uu$ is the velocity of the fluid, $p$ is its pressure and $\mu$ is its viscosity.

\begin{figure}[!h]
\begin{center}
\includegraphics[width=0.9\textwidth]{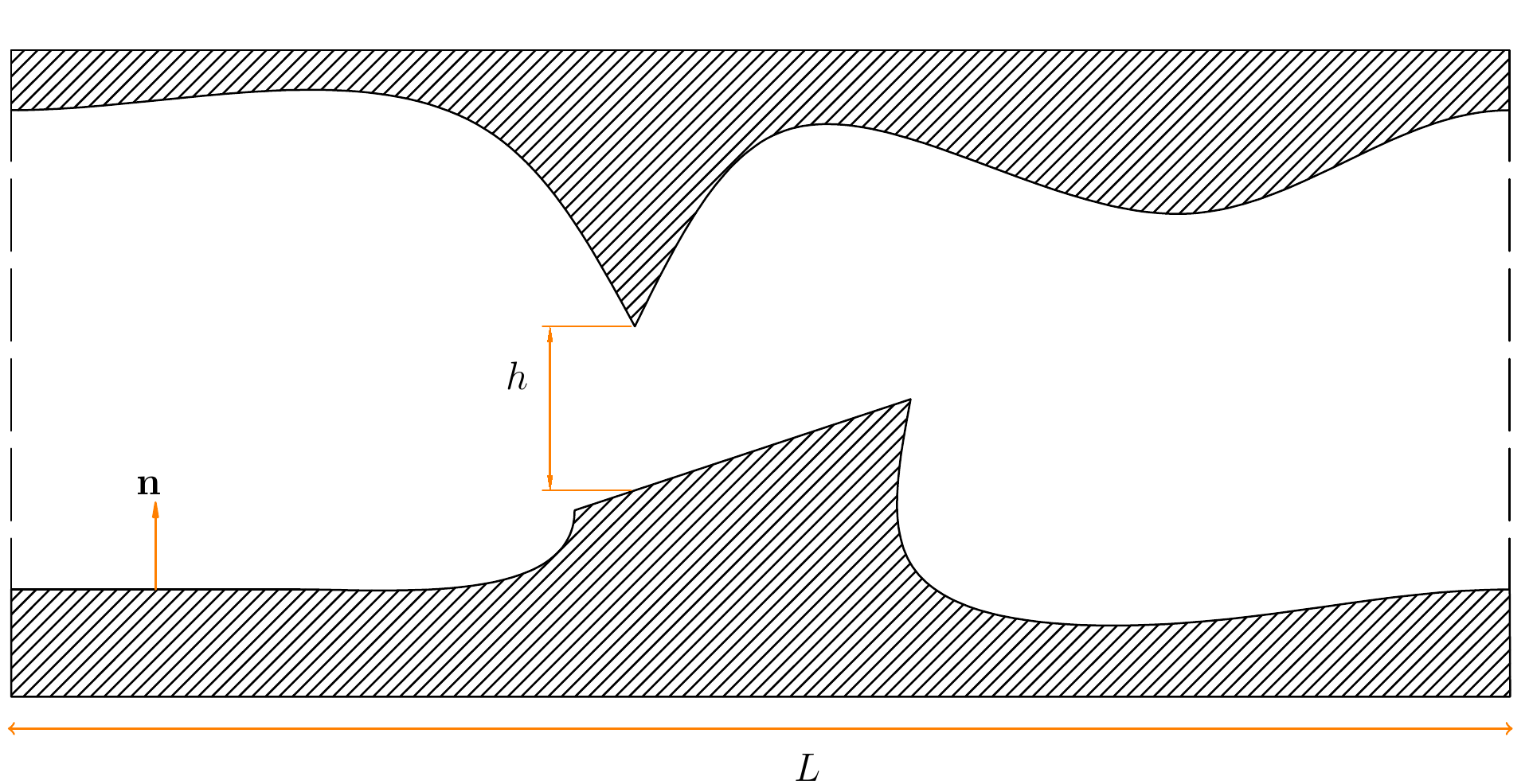}
\end{center}
\caption{Sketch of a periodic channel with a Lipschitz boundary. The period of the channel is $L$ in the $x_1$ direction, and the minimum height of the channel is $h$. The normal vector on the channel walls points into the bulk fluid.}\label{fig:channel_domain}
\end{figure}

In this paper we will restrict our attention to periodic channel flow as depicted in Figure \ref{fig:channel_domain}. For boundary conditions, we will prescribe Dirichlet conditions on the velocity, and enforce periodicity in the $x_1$ direction on the velocity. In addition, we impose a pressure drop across the reference cell, so that the pressure itself is not periodic, but the gradient of the pressure is,
\begin{subequations}\label{eq:stokes_periodic}
\begin{alignat}{2}
	\uu(\xx) &= \mathbf{g}(\xx),\qquad &\xx\in\Gamma,\\
	\uu(\xx) &= \uu(\xx + L), \qquad &\xx\in\Omega,\\
	p(\xx) - p(\xx + L) &= p_0 - p_1,&\qquad \xx\in\Omega.
\end{alignat}
\end{subequations}
To ensure the incompressibility condition is satisfied, the boundary data $\mathbf{g}$ must satisfy the compatibility condition,
\[ \int_\Gamma \uu\cdot\nn~\text{d}S = 0.\]

\section{Boundary Integral Formulation}\label{sec:bie}

For clarity of exposition, in this section and the next we will present the boundary integral formulation, and the treatment for the corners on a non-periodic domain. The periodicity will be addressed in Section \ref{sec:periodicity}.

\begin{figure}
\begin{center}
\includegraphics[width=0.3\textwidth]{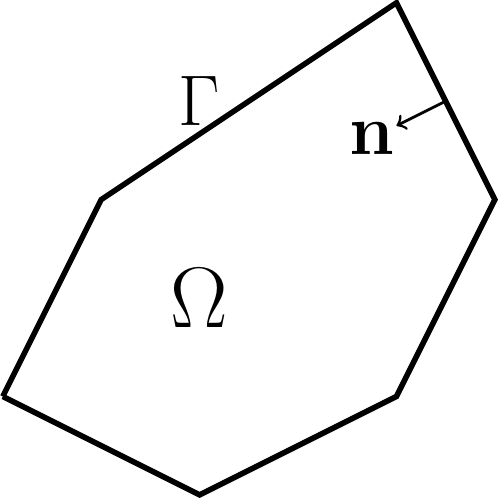}
\end{center}
\caption{Sketch of non-periodic, interior Lipschitz domain. The normal vector $\nn$ always points into the fluid.}\label{fig:simple_domain}
\end{figure}

Consider the Stokes equations \eqref{eq:stokes} defined inside a Lipschitz domain, along with Dirichlet boundary conditions on the velocity. For the non-periodic case, we will not be concerned with the pressure. The normal vector along the boundary $\Gamma$ always points into the fluid. See Figure \ref{fig:simple_domain} for a sketch of such a domain.

As given in \cite[Chapter 2]{Pozrikidis1992}, a solution to the forced Stokes equations
\[ -\mu\Delta\uu + \nabla p = \mathbf{g} \delta(\xx - \xx_0),\]
is given by the velocity and pressure pair
\[ u_j(\xx) = \frac{1}{4\pi\mu}G_{j\ell}(\xx - \xx_0)g_\ell, \qquad p(\xx) = \frac{1}{4\pi} p_j(\xx - \xx_0)g_j.\]
The stress tensor $\sigma_{j\ell}$ at a point $\xx$ is given by the combination of the velocity and pressure,
\begin{align*}
	 \sigma_{jm} &= \frac{1}{4\pi}\left(-\delta_{jm}p_\ell(\xx - \xx_0) + \frac{\partial G_{j\ell}}{\partial x_m}(\xx - \xx_0) +\frac{\partial G_{m\ell}}{\partial x_j}(\xx - \xx_0)\right)g_\ell\\
	 &= \frac{1}{4\pi} T_{j\ell m}(\xx - \xx_0)g_\ell.
\end{align*}
Here and in the remainder of this paper we have used the Einstein summation convention, where summation over repeated indices is implied. The letters $i$ and $k$ are reserved for later use and are therefore not used as indices.

In $\mathbb{R}^2$, we have that
\[ G_{j\ell}(\rr) = -\delta_{j\ell} \log(r) + \frac{r_j r_\ell}{r^2},\qquad p_j(\rr) = 2\frac{r_j}{r^2}, \qquad T_{j\ell m}(\rr) = -4\frac{r_j r_\ell r_m}{r^4}.\]
Here we have used $r$ to denote the Euclidean norm of $\rr := \xx  - \xx_0$.

The tensors $G_{j\ell}$ and $T_{j\ell m}$ are called the Stokeslet and stresslet respectively. From the Stokeslet and stresslet we can define the single- and double-layer potentials, $\mathbb{S}[\qq](\xx)$ and  $\mathbb{D}[\qq](\xx)$, as a convolution of the Stokeslet and stresslet with a density function $\qq(\xx)$ defined over $\Gamma$. Explicitly,
\begin{alignat*}{2}
	 \mathbb{S}[\qq](\xx) &= \frac{1}{4\pi\mu}\int_\Gamma q_j(\yy)G_{j \ell}(\xx - \yy)~\text{d}s(\yy),\qquad &\xx\in\Omega,\\
	 \mathbb{D}[\qq](\xx) &= -\frac{1}{4\pi}\int_\Gamma q_j(\yy)T_{j\ell m}(\xx - \yy) n_m(\yy)~\text{d}s(\yy),\qquad &\xx\in\Omega,
\end{alignat*}
where $\nn$ is the unit normal vector pointing into the fluid.

Following \cite{Hebeker1986}, we will express the solution to \eqref{eq:stokes} as a combination of the single- and double-layer potentials
\begin{equation}\label{eq:hebeker}\uu(\xx) =  2\mathbb{D}[\qq](\xx) + \eta \mathbb{S}[\qq](\xx) := \mathbb{K}[\qq](\xx), \qquad \xx\in\Omega,\end{equation}
where $\eta > 0$ is an arbitrary constant which governs the relation between the single and double layer potentials. To obtain a well-conditioned system $\eta$ cannot be too large; we will always take $\eta = 1$.

This equation is valid for $\xx\in\Omega$. To obtain a boundary integral equation (BIE) we will take the limit of \eqref{eq:hebeker} as $\xx$ approaches a point $\xx_0\in\Gamma$. To do this, we will need the limiting values of the single- and double-layer potentials,
\begin{align*}
	\lim\limits_{ \xx\to\xx_0\in\Gamma} \mathbb{S}[\qq](\xx) &= \mathbb{S}[\qq](\xx_0),\\
	\lim\limits_{ \xx\to\xx_0\in\Gamma} \mathbb{D}[\qq](\xx) &= -\frac{1}{2}\qq(\xx_0) + \mathbb{D}[\qq](\xx_0).
\end{align*}

Applying these limits to \eqref{eq:hebeker} and matching it to the boundary condition $\mathbf{g}$ yields the BIE
\begin{equation}\label{eq:bie} -\qq(\xx_0) + \mathbb{K}[\qq](\xx_0) = \mathbf{g}(\xx_0),\end{equation}
where $\xx_0\in\Gamma$, and $\mathbb{K}$ is the operator defined in \eqref{eq:hebeker} restricted to the boundary.

If $\Gamma$ is a Lyapunov curve (see \cite[Chapter 1]{Gohberg1992} for a precise definition), then  $\mathbb{K}$ is a compact operator, with eigenvalues accumulating at zero \cite[Chapter 15]{Kim2005}. In this case \eqref{eq:bie} can be analyzed using Fredholm theory. In particular the Fredholm alternative applies, and in \cite{Hebeker1986} this is used to demonstrate the existence and uniqueness of solutions. An alternative to the formulation \eqref{eq:hebeker} is the so-called completed double-layer formulation \cite{Power1993,Power1987}. This formulation has the advantage of not requiring the single-layer potential, which is weakly singular on the boundary. We have chosen \eqref{eq:hebeker} because numerical experiments \cite{Zinchenko2006} have shown it to be more stable for squeezing drops, and in the periodic case it allows us to naturally impose a global pressure gradient, as will be discussed in Section \ref{sec:periodicity}.

If, as in our case, $\Gamma$ is only Lipschitz smooth, then $\mathbb{K}$ is not compact, and Fredholm theory cannot be applied. Furthermore, the double-layer potential involves the normal vector, meaning that the double-layer kernel cannot be evaluated at any corner points on $\Gamma$. Nonetheless, it can be shown that \eqref{eq:bie} has a unique solution, even when $\Gamma$ is only Lipschitz continuous (see for example \cite{Verchota1984}). In this case the density function $\qq$ is singular at the corner points, however, \eqref{eq:hebeker} can still be evaluated for any $\xx \in \Omega$, and \eqref{eq:bie} can be evaluated anywhere on $\Gamma$ except at the corner points.

\section{Numerical Methods}

To evaluate \eqref{eq:bie}, we use the Nystr\"{o}m method, as described in \cite [Chapter 4]{Atkinson1997}. Let $\gamma(s),~s\in[0,2\pi]$ parameterize $\Gamma$. The BIE \eqref{eq:bie} can be written in the abstract form
\begin{equation}\label{eq:bie_abstract}
	-q_\ell(t) + \int_0^{2\pi} K_{j\ell}(t, s)q_j(s)|\gamma'(s)|~\text{d}s = g_\ell(t),  \qquad t \in [0,2\pi],
\end{equation}
where $K$ in this case is the kernel of $\mathbb{K}$, i,e.,
\begin{equation}\label{eq:kernel} K_{j\ell}(s,t) =  \frac{1}{2\pi}T_{j\ell m}(\gamma(t)-\gamma(s)) n_m(s)+ \frac{\eta}{4\pi\mu}G_{j\ell}(\gamma(t)-\gamma(s)).\end{equation}

We will approximate the integral in \eqref{eq:bie_abstract} using a composite Gauss-Legendre quadrature scheme of $n_\text{pan}$ panels and $n_q$ quadrature points per panel,
\begin{equation}\label{eq:discrete}
	-q_\ell(t) + \sum\limits_{n=1}^{N} K_{j\ell}(t,s_n)q_j(s_n)| \gamma'(s_n)| w_n = g_\ell(t), \qquad t\in[0,2\pi],
\end{equation}
where $N = n_\text{pan}n_q$ is the total number of quadrature points, and $w_n$ is the quadrature weight corresponding to the quadrature point $s_n$. For the remainder of this paper we will assume $n_q = 16$.  

We will then enforce \eqref{eq:discrete} at the quadrature points $s_m,~m=1,\hdots,N$ to get the linear system
\[ -q_\ell(s_m) + \sum\limits_{n=1}^{N} K_{j\ell}(s_m,s_n) q_j(s_n)|\gamma'(s_n)| w_n = g_\ell(s_m), \qquad m=1,\hdots, N.\]
Here the point values of the density function $q$ are unknown. When setting up the composite quadrature rule, care must be taken to ensure that each corner on $\Gamma$ is at the intersection of two panels. In this way we avoid difficulties arising from trying to evaluate the normal vector on a corner, as the Gauss-Legendre quadrature points cluster near but are never located at panel endpoints.

\subsection{Singular Quadrature}\label{sec:singular_quad}

When $t=s$, the kernel of the double-layer potential $T_{j\ell m}(0)n_m(s)$ has a removable singularity, provided that $t$ is not a corner point. The Stokeslet, however, must be handled using a specialized quadrature technique. We will use the approach given in \cite{Ojala2012}. We begin by writing the single-layer potential in complex variables,
\begin{equation}\label{eq:slp_complex}
\begin{aligned}\mathbb{S}[q](z) = \frac{1}{8\pi\mu}\int_\Gamma q(\tau)|\text{d}\tau| &+ \frac{1}{4\pi\mu}\int_\Gamma q(\tau)\log(|\tau - z|)|\text{d}\tau|\\& -\frac{1}{8\pi\mu}\int_\Gamma\overline{q(\tau)}\frac{\tau - z}{\overline{\tau} - \overline{z}}|\text{d}\tau|,\end{aligned}
\end{equation}
where we have introduced a slight abuse of notation to denote $q(z)$ as density function $\qq$ written as a complex number, i.~e., $z = x_1 + i x_2$ and $q(z) = q_1(\xx) + i q_2(\xx)$.

The only term in the complex formulation of the Stokeslet that does not have a finite limit as $\tau\to z$ is the term with the logarithm. We will consider this integral over a single panel, $\Gamma_k$. The panel extends from $\alpha_1$ to $\alpha_2$, with $\alpha_1, \alpha_2\in\mathbb{C}$, and is parameterized by $\alpha\in [\alpha_1, \alpha_2]$. Let $\alpha$ be parameterized as,
\[ \alpha(t) = \frac{\alpha_1 + \alpha_2}{2} + \psi t,\]
where $t \in [-1, 1]$ and $\psi = (\alpha_2 - \alpha_1)/2$. The $\log$ term in the integral above can be rewritten in the form,
\[ \log(|\tau(\alpha) - z|) = \log\left(\biggr| \frac{\psi(\tau(\alpha) - z)}{\alpha - \alpha_z}\frac{\alpha - \alpha_z}{\psi}\biggr|\right) = \log\left(\biggr|\psi\frac{\tau(\alpha) - z}{\alpha - \alpha_z}\biggr|\right) + \log\left(\biggr|\frac{\alpha - \alpha_z}{\psi}\biggr|\right),\]
where $\tau(\alpha_z) = z$. Define $t_z$ to be such that $\alpha(t_z) = \alpha_z$. Then we have that
\[ \frac{\alpha - \alpha_z}{\psi} = t - t_z,\]
which allows us to rewrite the $\log$ integral in \eqref{eq:slp_complex} as
\begin{align*} \int_{\Gamma_k} q(\tau)&\log(|\tau - z|)|\text{d}\tau| = \int_{\Gamma_k} q(\tau(\alpha))\log(|\tau(\alpha) - z|)|\tau'(\alpha)|~\text{d}\alpha\\
&=\int_{\Gamma_k} q(\alpha)\log\left(\biggr|\psi\frac{\tau(\alpha) - z}{\alpha - \alpha_z}\biggr|\right)~\text{d}\alpha + \int_{-1}^1 q(\alpha(t))\log(|t - t_z|)|\tau'(\alpha(t))|~\text{d}t.
\end{align*}
The first integral can be evaluated using the fact that
\[\lim\limits_{\alpha \to\alpha_z} \frac{\tau(\alpha) - z}{\alpha - \alpha_z} = \tau'(\alpha_z).\]

To evaluate the second integral, we expand the density $q(t)$ as a polynomial series to obtain,
\[I(t_z) = \int_{-1}^1 q(\alpha(t))\log(|t - t_z|)|\tau'(\alpha(t))|~\text{d}t \approx \sum\limits_{j=0}^{n_q - 1} c_j \int_{-1}^1  t^j  \log(|t - t_z|)|\tau'(\alpha(t))|~\text{d}t.\]

Denoting $\mathbf{h} = h_j, ~j=0,\hdots n_q - 1$ as
\[ h_j(t_z) = \int_{-1}^1 t^j  \log(|t - t_z|)|\tau'(\alpha(t))|~\text{d}t,\]
allows us to write $I$ in the compact form
\[ I(t_z) \approx \sum\limits_{j=0}^{n_q - 1} c_j h_j = \mathbf{c}^T \mathbf{h}(t_z).\]
The integrals in $\mathbf{h}$ can all be computed analytically using a recursion relation. The vector $\cc = \{c_j\}_{j=0}^{n_q - 1}$ can be computed by solving the Vandermonde system
\[ \bm{\vec{q}} = V \cc,\]
where $\bm{\vec{q}} = \{q(t_0), \hdots, q(t_{n_{q-1}})\}$, and $V$ is the Vandermonde matrix. Thus $I$ can be approximated by
\[ I(t_z) \approx (V^{-1} \bm{\vec{q}})^T\mathbf{h}(t_z) = \bm{\vec{q}}^T(V^{-1})^T \mathbf{h}(t_z) = \bm{\vec{q}}^T \bm{\omega}(t_z),\]
where $\bm{\omega}(t_z)$ is the solution to the linear system $V^T\bm{\omega}(t_z) = \hh(t_z)$. The stability of these computations is discussed in \cite [Appendix A]{Helsing2008b}, here we simply note that the computations are stable up to at least $n_q = 40$. Note that this linear system is independent of $\bm{\vec{q}}$. For now we will need the values of $I$ only at the quadrature points $t_0, \hdots, t_{n_{q-1}}$.  Since we have rescaled and rotated the panel $\Gamma_k$ to run from -1 to 1, the corrected weights $\bm{\omega}(t_z)$ can be precomputed for each of the quadrature points $t_0, \hdots, t_{n_{q-1}}$.

Quadrature points on adjacent panels will also need to be corrected. How many points to correct depends on accuracy considerations, but for $n_q = 16$, numerical results have shown that if all the panels are the same size (in parameter space), then correcting the four closest quadrature points in each adjacent panel is sufficient. Again, these corrections can be precomputed.

\subsection{Recursively Compressed Inverse Preconditioning}\label{sec:rcip}

\begin{figure}[!h]
\begin{center}
\begin{tabular}{c c c c}
\includegraphics[width=0.2\textwidth]{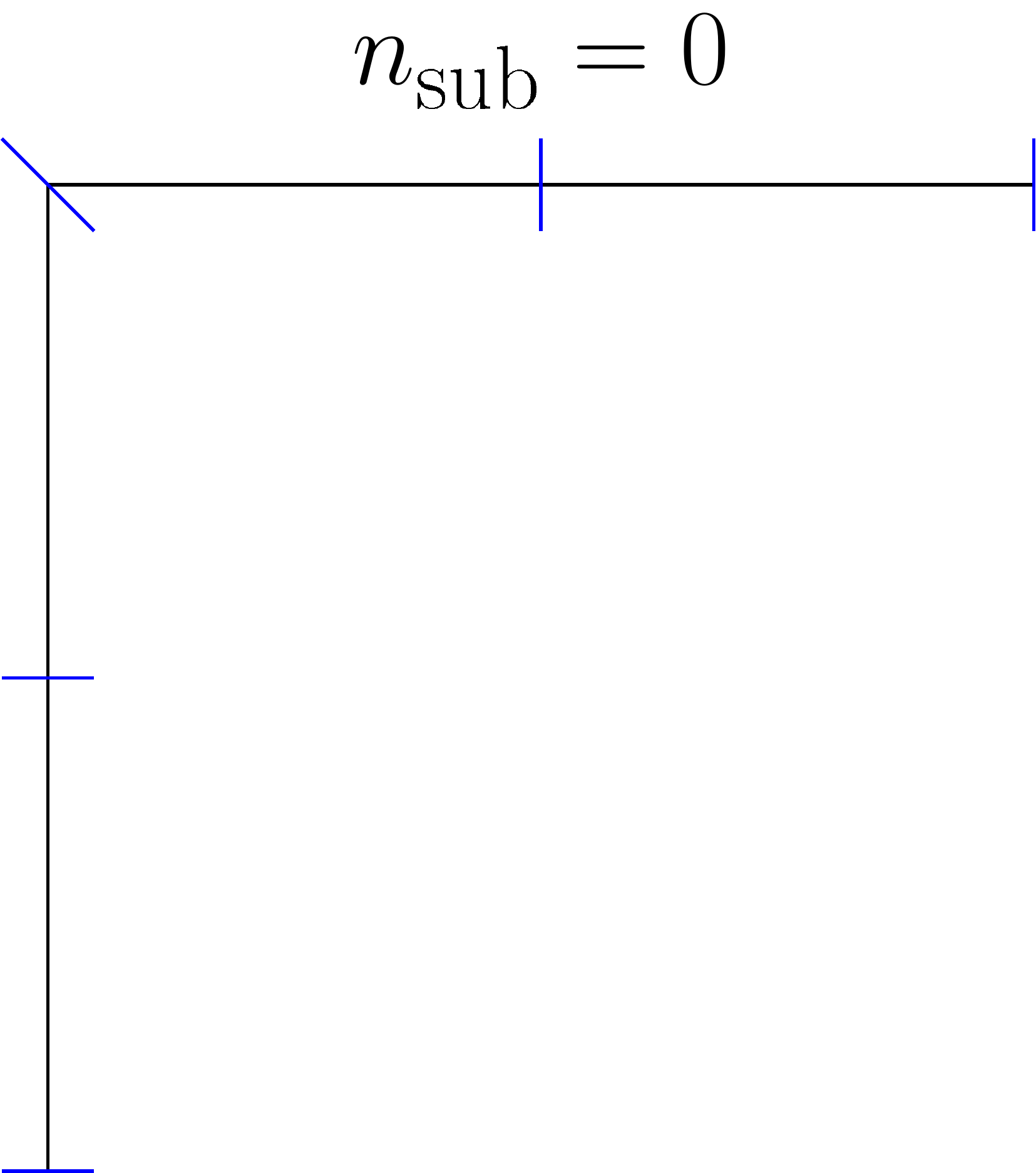} &
\includegraphics[width=0.2\textwidth]{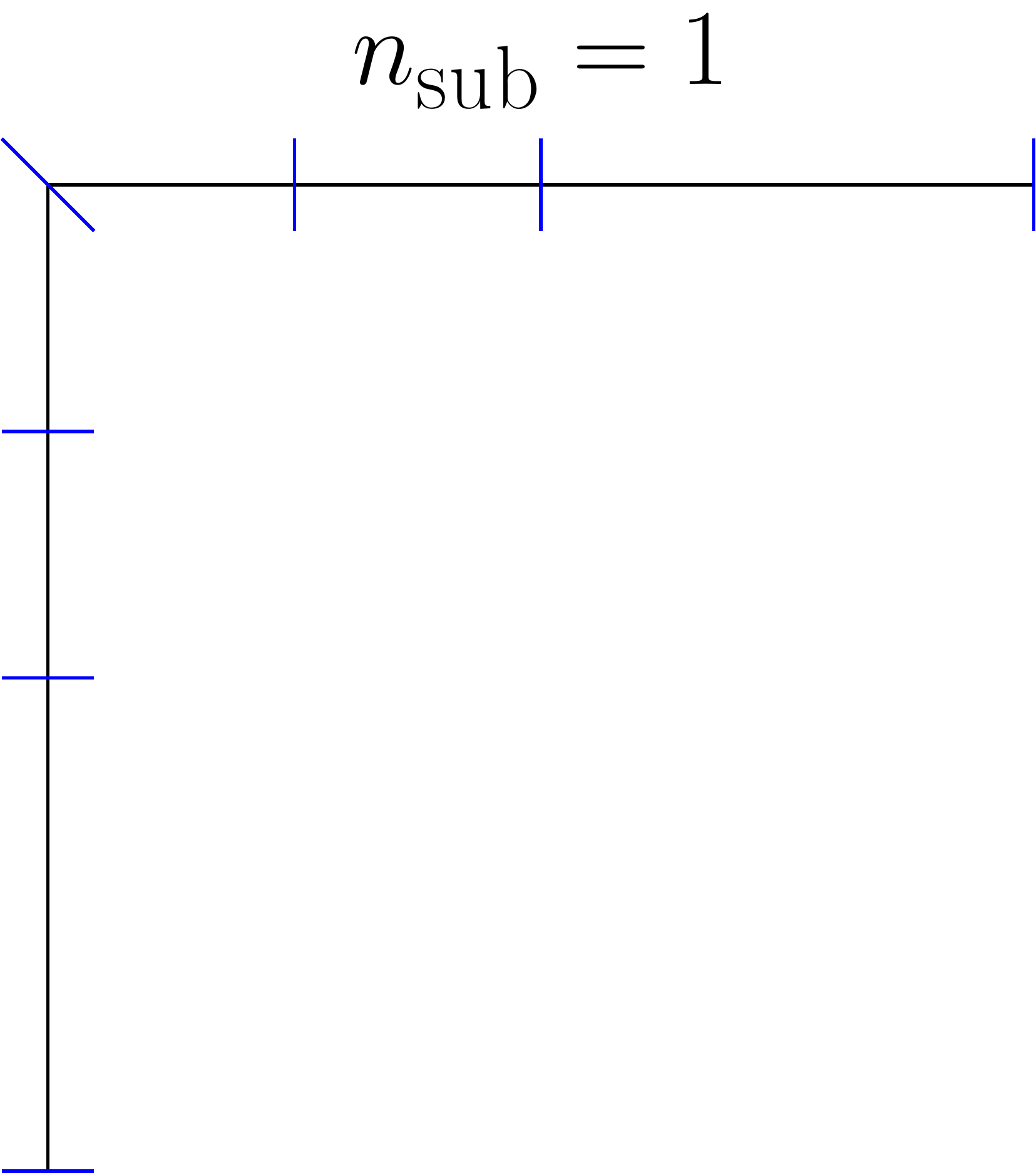} &
\includegraphics[width=0.2\textwidth]{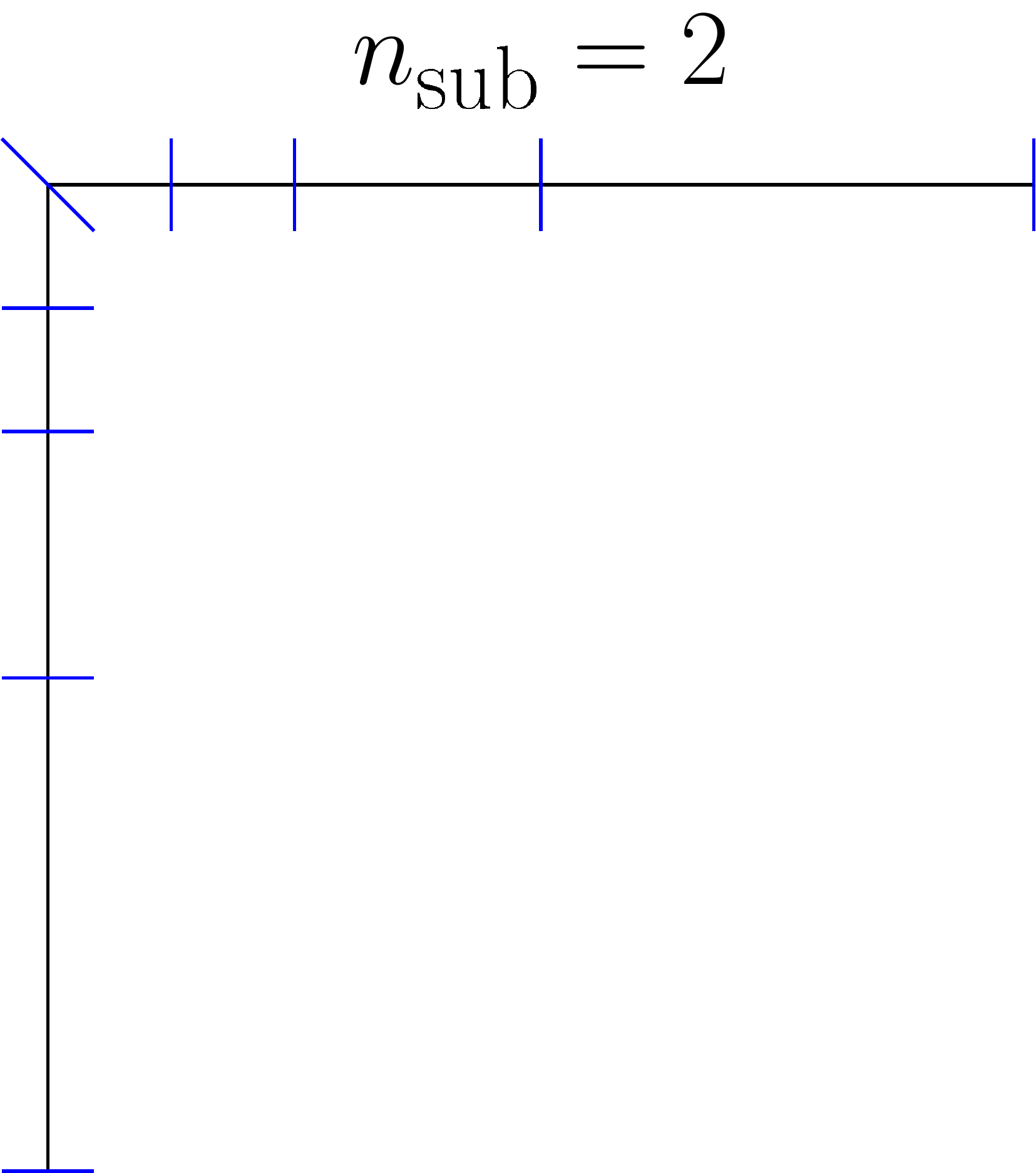} &
\includegraphics[width=0.2\textwidth]{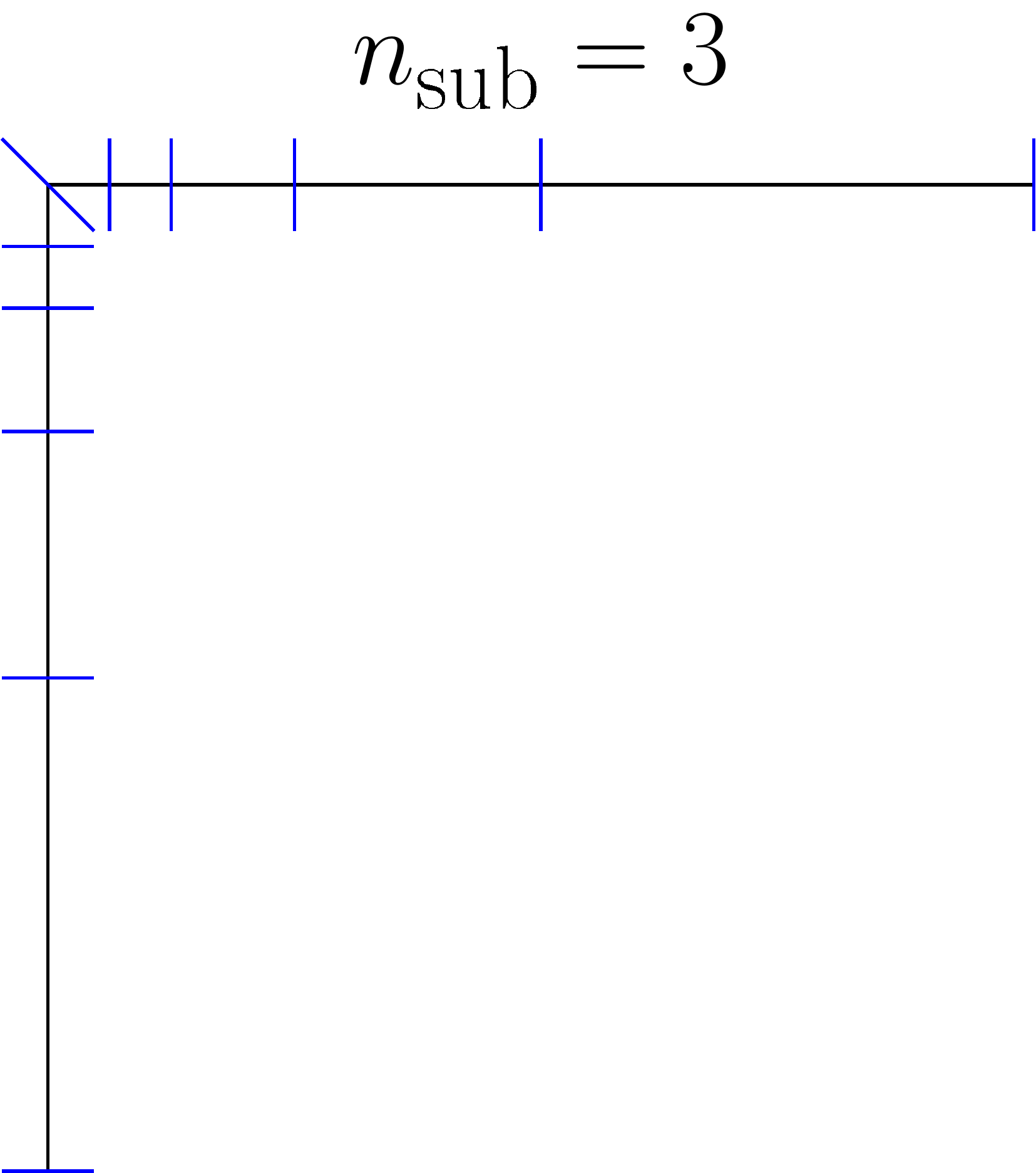}
\end{tabular}
\end{center}
\caption{Local panel refinement around a corner. Each addition subdivision adds two extra panels and $4n_q$ additional unknowns.}
\label{fig:localreffig}
\end{figure}

As previously mentioned, in the case of Lipschitz domains, \eqref{eq:bie} has a unique solution. However at the corners, the density function $\qq$ becomes singular. Therefore Gauss-Legendre quadrature fails to accurately integrate the integrands in \eqref{eq:bie} near the corners. Local panel refinement around the corners is one way to mitigate this issue; however, this approach adds a potentially large number of new unknowns, see Figure~\ref{fig:localreffig}. In addition, the condition number of the resulting locally refined linear system grows with the number of quadrature points.

An alternative approach, \emph{recursively compressed inverse preconditioning} (RCIP) \cite{Helsing2013b,Helsing2008a}, can achieve the same accuracy as local refinement by performing a simple precomputation. As the desired accuracy requirements increase, the precomputation grows linearly with $n_\sub$, however, both the size and the condition number of the final linear system remain fixed.

The main idea behind RCIP is to transform the density function $\qq$ in \eqref{eq:bie} into a transformed density $\tilde{\qq}$ that is piecewise smooth everywhere on $\Gamma$. Then Gauss-Legendre quadrature can be effectively used. To create this transformation, the operator $\mathbb{K}$ is written as
\begin{equation}\label{eq:operator_split}
	\mathbb{K} = \mathbb{K}^\circ + \mathbb{K}^*,
\end{equation}
 where $\mathbb{K}^\circ$ is compact away from the corners, and $\mathbb{K}^*$ describes the corner interactions (note that $\mathbb{K}^*$ is \emph{not} the adjoint of $\mathbb{K}$). These operators will be defined precisely shortly. We then define the transformed density as
\begin{equation}\label{eq:transformed_density}
	 \tilde{\qq} = (-\II + \mathbb{K}^*)\qq,
\end{equation}
where $\II$ is the identity operator.

Using the split \eqref{eq:operator_split}, and the transformed density \eqref{eq:transformed_density}, we can convert \eqref{eq:bie} to the transformed BIE for $\tilde{\qq}$,
\begin{equation}\label{eq:bie_transformed} (\II + \mathbb{K}^\circ\mathcal{R})\tilde{\qq}(\xx_0) = \mathbf{g}(\xx_0), \qquad \xx_0\in\Gamma,\end{equation}
where $\mathcal{R} = (-\II + \mathbb{K}^*)^{-1}$. If we assume $\mathbf{g}$ is piecewise smooth,  it can be immediately seen that $\tilde{\qq}$ must also be piecewise smooth. Since $\mathbb{K}^\circ$ is a smoothing operator, then $\mathbb{K}^\circ\mathcal{R}\tilde{\qq}$ will be smooth everywhere. Then, by contradiction, in order for \eqref{eq:bie_transformed} to hold, $\II\tilde{\qq} = \tilde{\qq}$ must be piecewise smooth.

\begin{figure}[!h]
\begin{tabular}{c c c}
\includegraphics[width=0.3\textwidth]{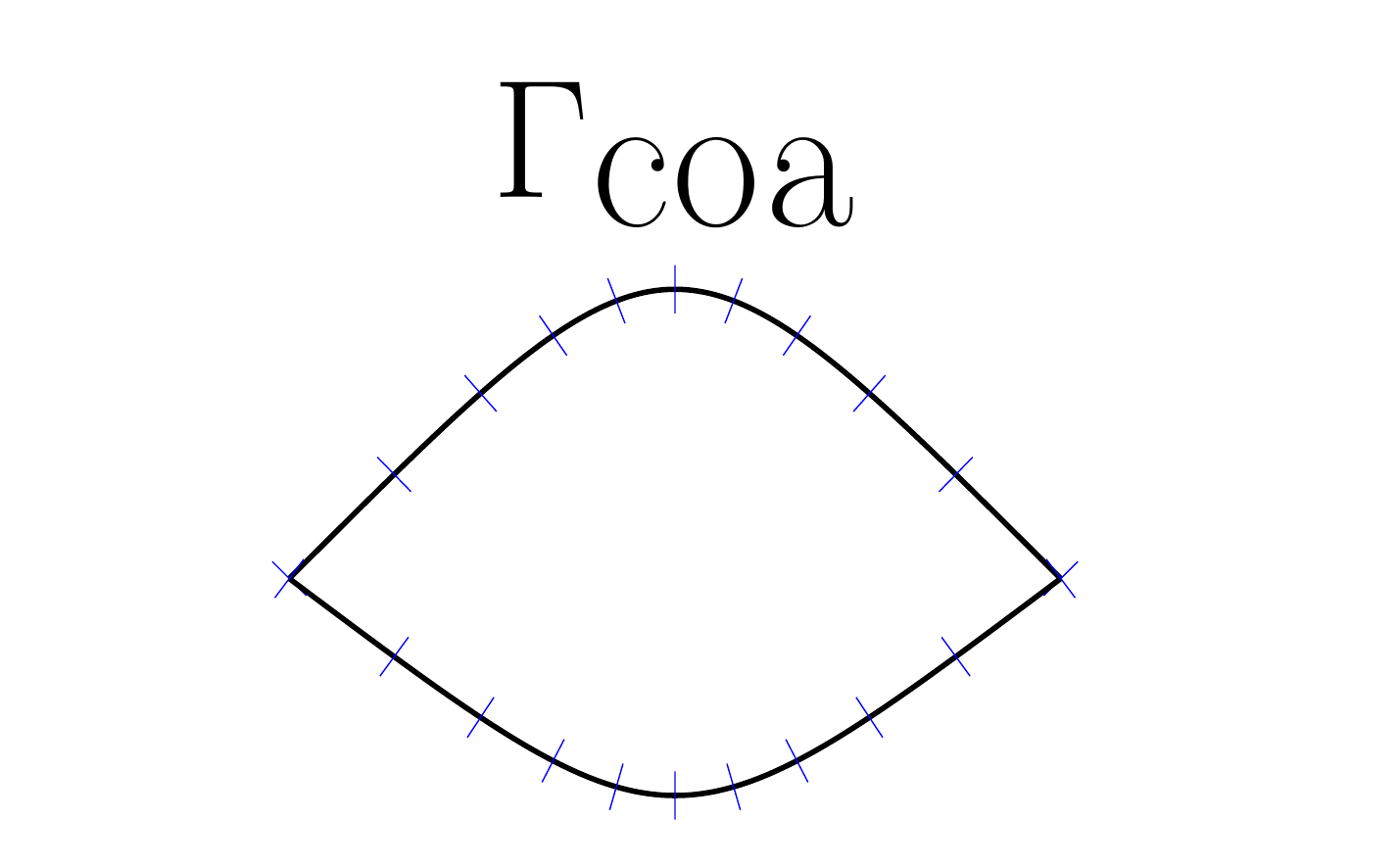} &
\includegraphics[width=0.3\textwidth]{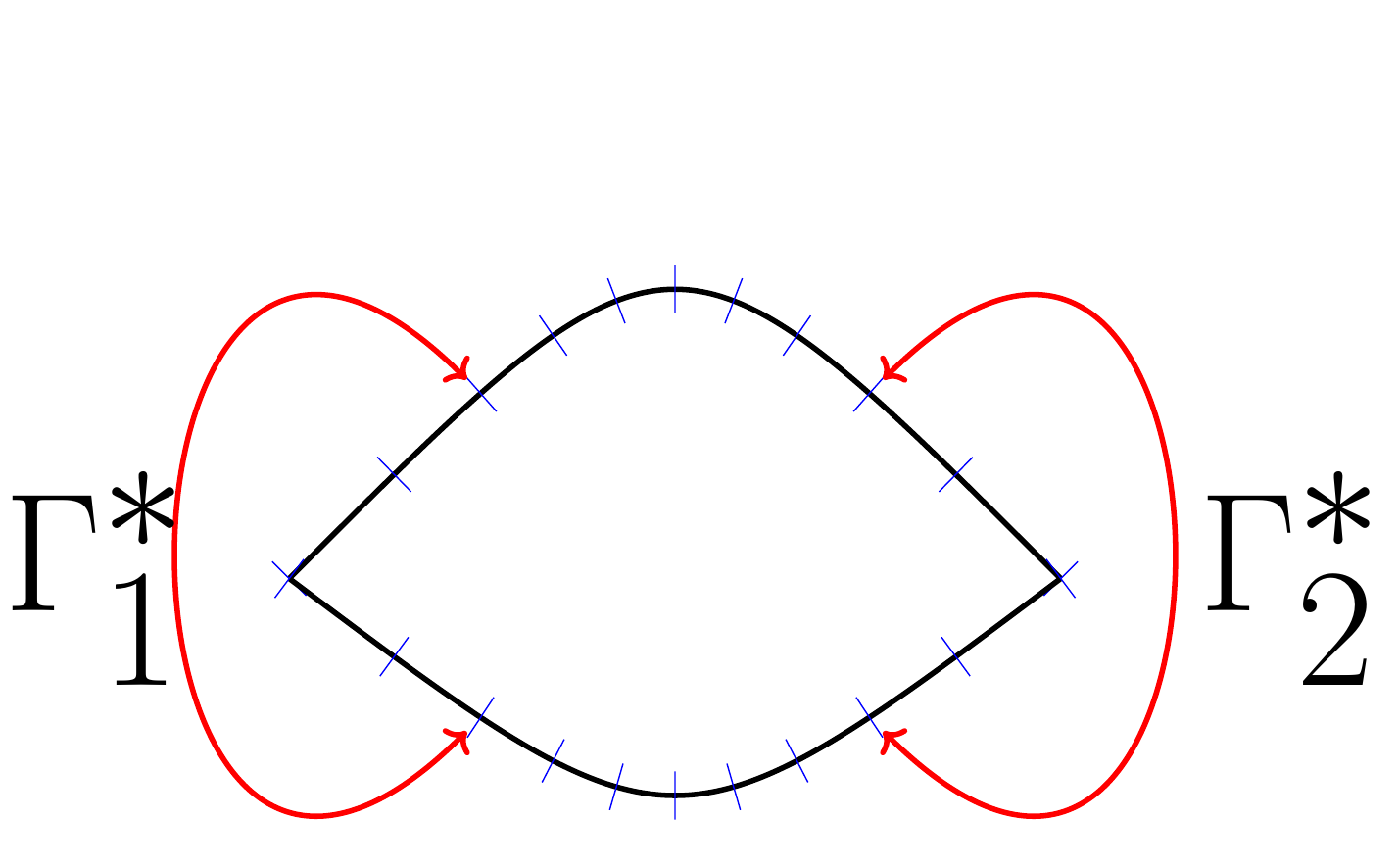} &
\includegraphics[width=0.3\textwidth]{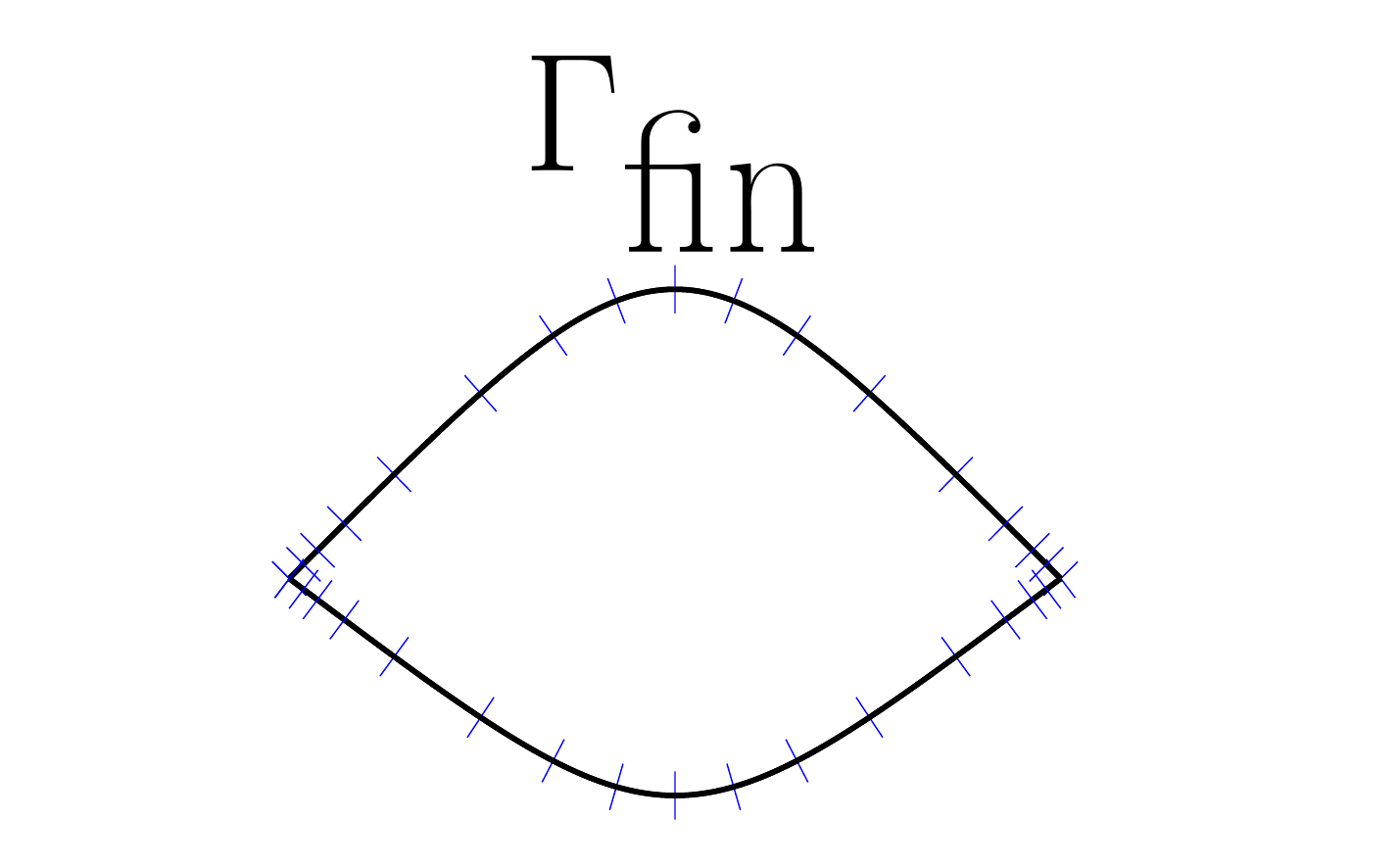}
\end{tabular}
\caption{Panel discretization using two meshes. The boundary $\Gamma$ is first discretized with a coarse mesh, $\Gamma_\coa$. The two panels on either side of each corner are denoted $\Gamma^*$, and the remaining panels are denoted $\Gamma^\circ$. The panels in $\Gamma^*$ closest to the corners are dyadically refined $n_\sub$ times to obtain $\Gamma_\fin$. Each panel, regardless of its size, contains $n_q$ Gauss-Legendre quadrature points.
}\label{fig:discretization}
\end{figure}

It remains to discretize \eqref{eq:bie_transformed}. One possibility for discretization is to use two meshes: a \emph{coarse} mesh $\Gamma_\coa$ on which $\mathbb{K}$ is discretized, and a \emph{fine} mesh $\Gamma_\fin$, on which $\mathcal{R}$ is discretized. To get $\Gamma_\fin$ from $\Gamma_\coa$ we first designate the two panels on either side of corner $j$ to be the subset $\Gamma^*_j\subset\Gamma_\coa$, and define $\Gamma^\circ$ to be $\Gamma_\coa \setminus \cup_{j=1}^{n_c} \Gamma^*_j$, where $n_c$ is the number of corners. The panels on either side of corner $j$ in $\Gamma^*_j$ are then dyadically refined $n_\sub$ times to get $\Gamma^*_{j,\fin}$. The fine mesh is defined as $\Gamma_\fin = \Gamma^\circ \cup\left(\cup_{j=1}^{n_c} \Gamma^*_{j,\fin}\right)$. An example of this is shown in Figure~\ref{fig:discretization}.

We will define $\mathbb{K}^\circ$ and $\mathbb{K}^*$ to be the operator $\mathbb{K}$ restricted to the domains $\Gamma^\circ$ and $\Gamma^*$ respectively. Note that $\mathbb{K}^\circ$ is a compact operator.  The operator $\mathbb{K}^\circ$ can be discretized on $\Gamma_\coa$ to get the matrix $\mathbf{K}^\circ$. This matrix is equivalent to the discretization of $\mathbb{K}$ on $\Gamma_\coa$, but with the entries corresponding to both source and target being outside $\Gamma^\circ$ set to zero.

\begin{figure}[!h]
\begin{center}
\includegraphics[width=0.8\textwidth]{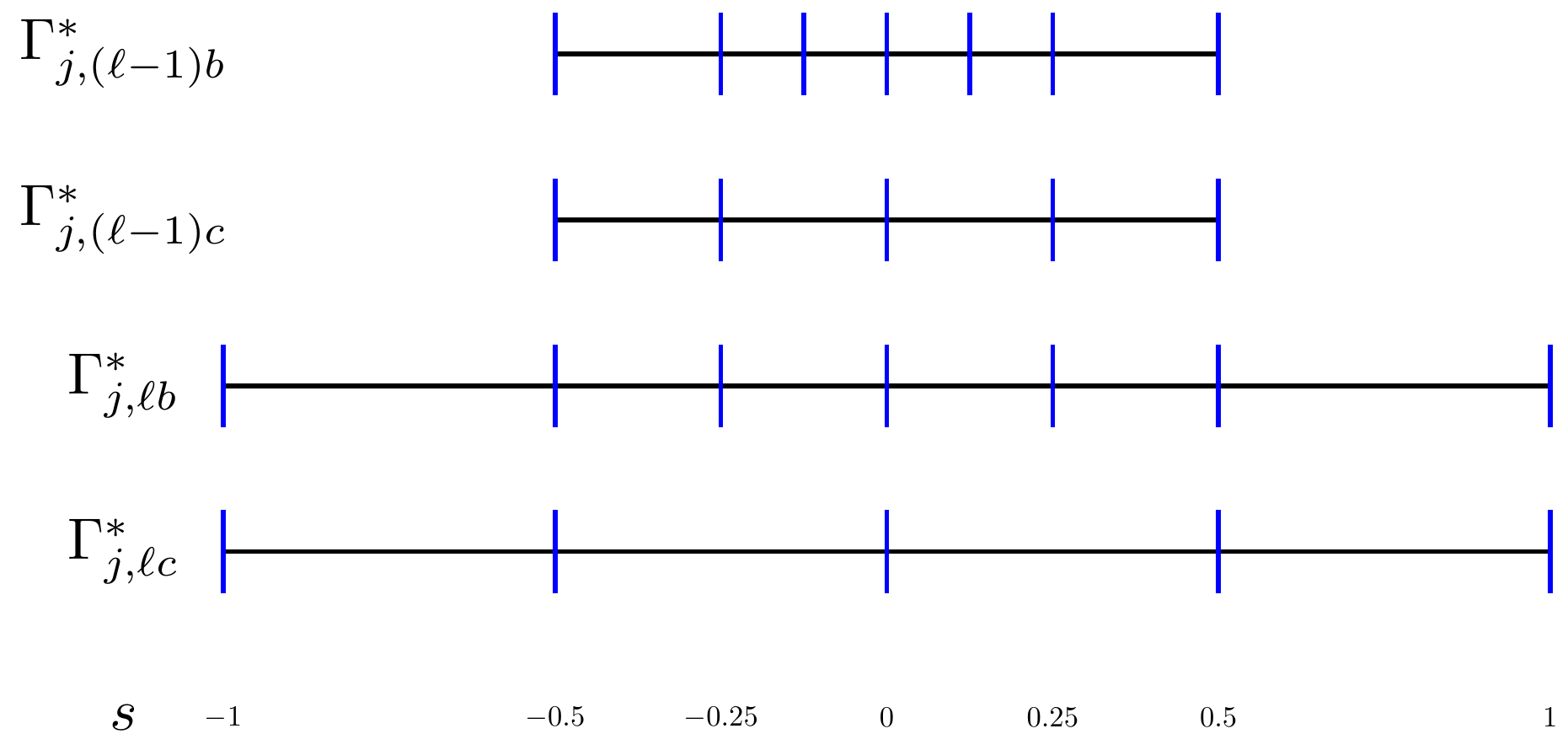}
\end{center}
\caption{Nested meshes in parameter space.  Let $\Gamma^*_{j,\ell}$ be parameterized by $s\in[-1,1]$, with corner $j$ at $s = 0$. Then the panel junctions on $\Gamma^*_{j,\ell c}$ will be at $s = \{-1, 0.5, 0, 0.5, 1\}$, and the panel junctions of $\Gamma^*_{j,\ell b}$ will be at $s = \{-1, -0.5, -0.25, 0, 0.25, 0.5, 1\}$. The panels are nested so that the junctions for $\Gamma^*_{j,(\ell-1)c}$ in the parameterization of $\Gamma^*_{j,\ell}$ will be at $s = \{-0.5, -0.25, 0, 0.25, 0.5\}$. }\label{fig:nested_meshes}
\end{figure}

The operator $\mathcal{R}$ could be discretized on $\Gamma_\fin$. This would, however, not be much use, since it would introduce a large number of unknowns. Instead we will exploit a forward recursion relation to construct $\RR$ on a sequence of meshes covering larger and larger portions of $\Gamma^*_j$. To define the recursion relation, we will need to provide some definitions of different types of meshes. For each corner $j$, define a sequence of meshes $\Gamma^*_{j,\ell},~\ell=1,\cdots, n_\sub$, with $\Gamma^*_{j,\ell-1}\subset\Gamma^*_{j,\ell}$ and $\Gamma^*_{j,n_\sub} = \Gamma^*_{j}$. On each $\Gamma^*_{j,\ell}$ we will have a six panel type-b mesh, $\Gamma^*_{j,\ell b}$, and a four panel type-c mesh $\Gamma^*_{j,\ell c}$. The type-b and type-c panels will be related as shown in Figure \ref{fig:nested_meshes}.

To interpolate between type-b and type-c meshes defined over the same interval, we introduce the prolongation matrix $\PP_{bc}$. In addition,  defining $\WW_b$ and $\WW_c$ to be a diagonal matrix containing quadrature weights on the type-b and type-c meshes respectively, we can define the weighted prolongation matrix $\PP_{Wbc}$ as
\[ \PP_{Wbc} = \WW_b \PP \WW_c^{-1}.\]

Let $\mathbf{K}_{j,\ell b}$ be the discretization of $\mathbb{K}$ on $\Gamma^*_{j,\ell b}$. This matrix will be of size $6n_q$ for scalar problems, or $12n_q$ for vector problems in $\mathbb{R}^2$. Define $\RR_{j,1}$ as
\[ \RR_{j,1} = \PP_{Wbc}^T\left(-\ii_b + \mathbf{K}_{j,1b}\right)\PP_{bc}.\]

We then compute the sequence $\RR_{j,2}, \hdots, \RR_{j,n_\sub}$ using the recursion relation
\begin{equation}\label{eq:recursion} \RR_{j,\ell} = \PP_{Wbc}^T\left(\mathbb{F}\{\RR_{j,\ell-1}^{-1}\} - \ii_b^\circ + \mathbf{K}^\circ_{j,\ell b}\right)^{-1}\PP_{bc},\end{equation}
where $\ii_b^\circ$ and $\mathbf{K}^\circ_{j,\ell b}$ are the identity matrix and $\mathbf{K}_{j,\ell b}$ with the entries in the two panels around the corner zeroed out respectively. The operator $\mathbb{F}\{\cdot\}$ zero pads the argument to turn a matrix defined on $\Gamma^*_{j,\ell c}$ into one defined on $\Gamma^*_{j,(\ell +1)b}$.

Once we have the matrices $\mathbf{R}_{j,n_\sub}$, we can construct $\hat{\RR}$, the discretization of $\mathcal{R}$ on $\Gamma_\coa$. Outside of $\Gamma^*$, $\mathbb{K}_*$ is zero, so from the definition of $\mathcal{R}$, we obtain that $\hat{\RR}$ must be the identity matrix over $\Gamma^\circ$. Over $\Gamma^*_j$, the discretization of $\mathcal{R}$ is just $\mathbf{R}_{j,n_\sub}$.

To handle the log singularity when assembling the matrix $\mathbf{K}_{j,1b}$ and $\mathbf{K}^\circ_{j,\ell b}$, the same techniques as described in Section \ref{sec:singular_quad} can be used. Some bookkeeping is needed to account for the fact that the panels are not equally sized in parameter space. In some cases adjacent panels will be double or half the size; however, the modifications for the quadrature weights can still be precomputed for each case.

For wedge shaped corners, the domain segments $\Gamma_{j,\ell b}$ are self similar for $\ell = 1, \hdots, n_\sub$ and $j = 1, \hdots, n_c$. If the operator $\mathbb{K}$ is scale invariant, then $\mathbf{K}^\circ_{j,\ell b}$ will be independent of $\ell$, and the recursion relation \eqref{eq:recursion} becomes a fixed point iteration which can be iterated to find $\RR_j$ to a desired tolerance without specifying $n_\sub$ \cite[Section  12]{Helsing2013b} . Unfortunately, in our case, because the single-layer potential contains a logarithmic term, $\mathbb{K}$ is not scale invariant, so we cannot use this idea, and $n_\sub$ must be specified a priori. A formulation involving just the double-layer potential \cite{Power1993,Power1987} would not suffer from this drawback and could be formulated as fixed point iteration. 

\subsection{Fast multiplication}\label{sec:multiplication}

The fully discrete transformed BIE defined on $\Gamma_\coa$ is
\begin{equation}\label{eq:bie_discrete_transformed}
	\left(\ii + \mathbf{K}^\circ\hat{\RR}\right)\tilde{\qq} = \mathbf{g},
\end{equation}
which can be solved using GMRES. To accelerate the required matrix-vector products from $\mathcal{O}(N^2)$ to something computationally feasible, fast summation methods are a necessity. To do this, we will rewrite $\mathbf{K}$ and $\hat{\RR}$ as the block matrices
\[ \mathbf{K} = \begin{pmatrix} \mathbf{K}^{**} & \mathbf{K}^{*\circ}\\ \mathbf{K}^{\circ*} & \mathbf{K}^{\circ\circ}\end{pmatrix}, \qquad \hat{\RR} = \begin{pmatrix} \RR & \mathbf{0} \\ \mathbf{0} & \ii^\circ\end{pmatrix}.\]

It follows that $\mathbf{K}^\circ$ can be expressed as
\[ \mathbf{K}^\circ = \mathbf{K} - \begin{pmatrix} \mathbf{K}^{**} & \mathbf{0} \\ \mathbf{0} &\mathbf{0}\end{pmatrix} = \mathbf{K} - \mathbf{K}^*.\]

Then \eqref{eq:bie_discrete_transformed} can be rewritten as
\[ \left( \ii + \mathbf{K}\hat{\RR} - \mathbf{K}^*\hat{\RR}\right)\tilde{\qq} = \mathbf{g}.\]

The matrix vector product $\hat{\RR}\tilde{\qq}$ can be done directly, since $\hat{\RR}$ is just the identity matrix with one $8n_q \times 8n_q$ block for each corner. The remaining matrix vector products $\mathbf{K}(\hat{\RR}\tilde{\qq})$ and $\mathbf{K}^*(\hat{\RR}\tilde{\qq})$ can be computed with a fast summation method, reducing the cost from $\mathcal{O}(N^2)$ to $\mathcal{O}(N\log N)$ or $\mathcal{O}(N)$, depending on the choice of fast summation method. We will use an $\mathcal{O}(N\log N)$ method that can naturally handle periodic problems, as  will be described in Section \ref{sec:periodicity}. 

\subsection{Post-processing} \label{sec:post}

After we have computed $\tilde{\qq}$, we evaluate the velocity at a point $\uu\in\Omega$ by using the regular quadrature rule
\[ u_\ell(\xx) \approx \sum\limits_{n=1}^N K_{j\ell}(\xx, s_n)\hat{q}_j(s_n)|\gamma'(s_n)|w_n,\]
where $\hat{\qq}=\hat{\RR}\tilde{\qq}$, and $K_{j\ell}$ is defined in \eqref{eq:kernel}. This quadrature is as accurate as if we had access to the fine density defined on $\Gamma_\fin$. From the definitions of $\tilde{\qq}$ and $\hat{\RR}$, it is clear that outside of $\Gamma^*$, $\hat{\qq}$ and $\qq$ are equivalent. On $\Gamma^*$, the transformed density $\hat{\qq}$ can be thought of as a weight corrected density function.

Both the double- and single-layer potentials contain integrals that become nearly singular when $|t-s|$ is small. This causes the numerical errors from standard Gauss-Legendre quadrature to grow quite large when evaluating the integrals at any point close to a boundary. To handle this, the integrals are treated in the manner described in \cite{Palsson2019b}. In short, the main idea is similar to that explained in Section~\ref{sec:singular_quad} where the density is expanded as a polynomial series, and the integrals are computed analytically using recursion relations.

\section{Periodicity}\label{sec:periodicity}

We now address the periodicity in \eqref{eq:stokes_periodic}. The periodic single- and double-layer potentials, $\mathbb{S}^P[\qq](\xx)$ and $\mathbb{D}^P[\qq](\xx)$ are defined as an infinite sum of the single- and double-layer potentials defined in Section \ref{sec:bie},
\begin{subequations}\label{eq:bie_infinite}
\begin{alignat}{2}
	\mathbb{S}^P[\qq](\xx) &= \frac{1}{4\pi\mu}\sum\limits_{\pp\in\mathbb{Z}^2} \int_\Gamma q_j(\yy) G_{j\ell}(\xx - \yy - L\pp)~\text{d}s(\yy), &\qquad \xx\in\Omega,\\
	\mathbb{D}^P[\qq](\xx) &= \frac{1}{4\pi}\sum\limits_{\pp\in\mathbb{Z}^2} \int_\Gamma q_j(\yy) T_{j\ell m}(\xx - \yy - L\pp)n_m(\yy)~\text{d}s(\yy), &\qquad \xx\in\Omega.
\end{alignat}
\end{subequations}

The interpretation of these periodic sums will be given in Section \ref{sec:ewald}. With these operators we can define the periodic operator $\mathbb{K}^P[\qq]$,
\begin{equation}\label{eq:beta_p}
\mathbb{K}^P[\qq](\xx) = 2\mathbb{D}^P[\qq](\xx) + \eta\mathbb{S}^P[\qq](\xx).
\end{equation}

Note that periodic sum is over $\pp\in\mathbb{Z}^2$; in other words we are implying periodicity in both spatial dimensions, even though in our actual problems the periodicity will be only in one direction. The reason for this is that it allows us to exploit a more efficient fast-summation method. To enforce the periodicity in one dimension we will embed the channel, which is only periodic over a length $L_1$ in the $x_1$ direction, in a doubly periodic box of size $L_1\times L_2$. A similar approach is used in \cite{Zhao2010} where a periodic flow in one direction is embedded in a two-dimensional periodic box. For simplicity we will assume here that $L_1 = L_2 = L$, but this need not be the case, and only minor modifications are needed to use a rectangular periodic box.

As mentioned in Section \ref{sec:gov}, a constant pressure drop $p_0 - p_1$ in the $x_1$ direction is applied. In order to impose this we will add on an unknown mean velocity $\langle \uu \rangle$ to the layer formulation \eqref{eq:bie},
\[ \uu(\xx) = \mathbb{K}^P[\qq](\xx) + \langle \uu \rangle.\]
Note that  $\langle \cdot\rangle$ denotes a volume average over a full periodic cell, including the parts outside $\Omega$. This quantity has no physical meaning beyond imposing a pressure gradient. For a channel with flat parallel walls, a simple alternative would be to add a predetermined background flow, as done in \cite{Palsson2019b}. 

The mean pressure gradient $\langle \nabla p\rangle$ must balance the net effect of the wall friction \cite{Zhao2010}.  From \eqref{eq:beta_p} the wall friction force is equal to $\eta \int_\Gamma \qq~\text{d}\Gamma$ \cite{Hebeker1986}. The system we have to solve is thus given by
\begin{align*}
	-\qq(\xx) + \mathbb{K}^P[\qq](\xx) + \langle \uu \rangle &= \mathbf{g}(\xx), \qquad \xx \in \Gamma,\\
	-\frac{1}{L^2}\eta \int_\Gamma q_1 ~\text{d}\Gamma &= \langle\nabla p\rangle, \\
	\int_\Gamma q_2 ~\text{d}\Gamma &= 0.
\end{align*}

The splits for the two dimensional Stokeslet and stresslet, as well as truncation estimates are derived in [26]. Here we list the results.

To compute the infinite sums in \eqref{eq:bie_infinite} we will use the spectral Ewald  method \cite{Lindbo2010,Lindbo2011}.  In the spectral Ewald method the infinite sums in $\mathbb{S}^P$ and $\mathbb{D}^P$  are split into two parts based on a so-called Ewald decomposition: one that decays exponentially fast in real space, and one that decays exponentially fast in Fourier space. The real space sum can be truncated in real space, while the Fourier sum can be truncated in Fourier space. The splits for the two-dimensional Stokeslet and stresslet, as well as truncation estimates are derived in \cite{Palsson2019b}. Here we list the results.

\subsection{Ewald splits}\label{sec:ewald}

The discretized periodic single- and double layer potentials are given by
\begin{align*}
	 \mathbb{S}^P[\qq](\xx) &\approx u^G_\ell(\xx) = \frac{1}{4\pi\mu}\sum\limits_{\pp\in\mathbb{Z}^2}^*\sum\limits_{m=1}^M  G_{j\ell}(\xx - \xx_m -L\pp)q_j(\xx_m)w_m,\\
	 \mathbb{D}^P[\qq](\xx) &\approx u^T_\ell(\xx) = \frac{1}{4\pi}\sum\limits_{\pp\in\mathbb{Z}^2}^*\sum\limits_{m=1}^M T_{j\ell s}(\xx - \xx_m - L \pp)q_j(\xx_m)n_s(\xx_m)w_m,
\end{align*}
where the * in the summations over $\pp$ indicates that we are excluding any $\pp$ that makes $\xx - \xx_m - L\pp = \mathbf{0}$. As written, these sums are not well-defined, as they are not convergent. By assuming that a pressure gradient is balancing the forces acting on the fluid, well-defined but slowly converging sums can be found in Fourier space. See \cite[Section 4]{Palsson2019b} for a discussion. Below we will introduce an Ewald decomposition for these sums that yields the same well-defined results, but achieves them by a split into two rapidly converging sums: one in real space, and one in Fourier space.

The discretized single-layer potential can be rewritten as the split
\begin{equation*}\begin{aligned}
	u^G_\ell(\xx) = \sum\limits_{\pp\in\mathbb{Z}^2}^* \sum\limits_{m=1}^M &G_{j\ell}^R(\xx - \xx_m + L\pp,\xi)q_j(\xx_m)w_m \\&+ \frac{1}{L^2}\sum\limits_{\kk\ne 0}\hat{G}^F_{j\ell}(\kk, \xi)\sum\limits_{m=1}^M e^{-i\kk\cdot(\xx - \xx_m)}q_j(\xx_m) w_m,
	\end{aligned}
\end{equation*}
where
\begin{align*}
G_{j\ell}^R(\rr,\xi)&= \frac{1}{4\pi}\left(e^{-\xi^2 r^2}\left(\frac{r_j r_\ell}{r^2}-\delta_{j\ell}\right) + \frac{\delta_{j\ell}}{2}E_1(\xi^2 r^2)\right)\\
\hat{G}_{j\ell}^F(\kk,\xi) &= \frac{1}{k^2}e^{-k^2 / (4\xi^2)}\left(\delta_{j\ell} - \frac{k_j k_\ell}{k^2}\right)\left(1 + \frac{k^2}{4\xi^2}\right), \qquad \kk\ne \mathbf{0}, \; k=|\mathbf{k}|.
\end{align*}
The $\kk=\mathbf{0}$ mode in the Fourier expansion has been shown to be 0 \cite{Pozrikidis1996}. The decomposition parameter $\xi > 0$ is called the Ewald parameter and determines the relative sizes of the real and Fourier parts of $\uu^G(\xx)$. Note that $\uu^G(\xx)$ itself is independent of $\xi$.

When applying the Nystr\"{o}m method to solve for the unknown pointwise density values, we do not wish to evaluate singular cases where $\xx - \xx_n + L\pp = \mathbf{0}$. This contribution to $\uu^G(\xx)$ can be skipped in the real space sum, but for the Fourier sum we will have to subtract off the limiting value given by
\[ \lim\limits_{r\to 0} (G_{j\ell}(\rr) - G^R_{j\ell}(\rr))q_j(\xx_n)w_n = \left(\frac{1}{2}\gamma + \log(\xi) + 1\right)q_\ell(\xx_n) w_n,\]
where $\gamma$ is the Euler-Mascheroni constant.

For the discretized double-layer potential, we use the split,
\begin{equation*}\begin{aligned}
	u^T_\ell(\xx) = \sum\limits_{\pp\in\mathbb{Z}^2}^* \sum\limits_{m=1}^M &T_{j\ell s}^R(\xx - \xx_m + L\pp,\xi)q_j(\xx_m)n_s(\xx_m)w_m \\&+ \frac{1}{L^2}\sum\limits_{\kk\ne 0}\hat{T}^F_{j\ell s}(\kk, \xi)\sum\limits_{m=1}^M e^{-i\kk\cdot(\xx - \xx_m)}q_j(\xx_m) n_s(\xx_m) w_m \\
	& + \frac{1}{L^2}\sum\limits_{m=1}^M \hat{T}^{F,0}_{j\ell s}(\xx_m) q_j(\xx_m)n_s(\xx_m) w_m.
	\end{aligned}
\end{equation*}
where
\begin{align*}
	T_{j\ell s}^R(\rr,\xi) &= \frac{1}{4\pi}e^{-\xi^2 r^2}\left( 2 \xi^2 (\delta_{j\ell} r_s + \delta_{js}r_\ell + \delta_{\ell s}r_j) - \frac{4 r_j r_\ell r_s}{r^4}(1 + \xi^2 r^2)\right),\\
	\hat{T}^F_{j\ell s}(\kk, \xi) &= \frac{i}{k^2} e^{-k^2/(4\xi^2)}\left(\delta_{j\ell}k_s + \delta_{j s}k_\ell + \delta_{\ell s}k_j) - \frac{2 k_j k_\ell k_s}{k^2}\right)\left(1 + \frac{k^2}{4\xi^2}\right), \qquad \kk \ne \mathbf{0}.
\end{align*}

For the $\kk = \mathbf{0}$ mode, the choice
\begin{equation}\label{eq:stresslet_zero_mode} \hat{T}^{F,0}_{j\ell s}(\xx) = \delta_{\ell s}x_j,\end{equation}
guarantees zero mean-flow though the reference cell \cite{Klinteberg2017}.

For the stresslet, the limit
\[ \lim_{r\to 0} (T_{j\ell s}(\rr) - T^R_{j\ell s}(\rr)) q_\ell n_s = \mathbf{0},\]
so no limiting value needs to be subtracted when $\rr = \mathbf{0}$.

To compute these sums, it is necessary to truncate them. Fortunately, both these sums decay exponentially fast, and they can be truncated to a desired tolerance following the estimates in \cite{Palsson2019b}. With appropriate scaling of the decomposition parameter, the real space part is computed in $\mathcal{O}(N)$ time and the Fourier space part is accelerated to $\mathcal{O}(N\log N)$ using the fast Fourier transform. In order to use FFTs, the source points are spread to a uniform grid where the computations for the Fourier space sum are carried out. The result is then gathered from the uniform grid to the target points. The spreading and gathering is done using truncated Gaussians whose shape parameter is selected to minimize the approximation error for a given support. For efficiency, fast Gaussian gridding \cite{Greengard2004} can be used in both the spreading and the gathering steps. For more details see \cite{Palsson2019b}.

A numerical difficulty in the periodic formulation lies in the fact that the Ewald representation of the stresslet is not translation invariant due to the zero mode \eqref{eq:stresslet_zero_mode}. This means that the submatrices $\RR_j$ are different for each corner, even if the corners have the exact same shape. Therefore to assemble $\mathbf{R}_j$ we first assemble the translation invariant part, and then add on \eqref{eq:stresslet_zero_mode} only at the end. Additionally, to avoid roundoff error, for large $n_\sub$, (where the smallest $\Gamma^*_{j,\ell}$ is $\mathcal{O}(10^{-16})$ or less) we introduce a local coordinate system for each corner centered at $\xx = \mathbf{0}$.

\section{Numerical Examples}

To test the periodization scheme, we can compare to a known exact solution. Pressure driven pipe flow creates the well-known Poiseuille parabolic flow profile in the $x_2$ direction. The exact solution for  flow through a flat channel with top wall at $x_2 = 1$ and bottom wall at $x_2 = 0$ is given by
\[ \mathbf{u}^{\text{exact}}(\xx) = \langle -\frac{(p_1 - p_0) x_2(x_2 - 1) }{2 L\mu}, 0 \rangle.\]
Note that this flow is constant in the $x_1$ direction, and therefore also periodic in the $x_1$ direction. As can be seen in Table \ref{tab:pipe_flow}, using only 8 panels per wall allows us to achieve a relative error of $10^{-11}$ everywhere up to a distance of $10^{-3}$ from the boundary. 

\begin{table}[!h]
\begin{center}
	\begin{tabular}{c | c  c c c }
	$n_{\text{pan}}$ per wall  & \multicolumn{4}{ c}{Distance from wall}\\
	\hline
	 & 0.5 & 0.1 & 0.01 & 0.001\\
	\hline
	1 & $1.64\times 10^{-1}$ & $2.07\times 10^{0}$ & $2.10\times 10^{1}$ & $8.44\times 10^{1}$ \\
	2 & $1.14\times 10^{-4}$ & $8.15\times 10^{-1}$ & $4.65\times 10^{0}$& $1.28\times 10^{1}$ \\
	4 & $5.95\times 10^{-10}$ & $1.12\times 10^{-9}$ & $5.03\times 10^{-9}$ & $ 7.07\times 10^{-9}$ \\
	8 & $8.88\times 10^{-15}$ & $1.25\times 10^{-13}$ & $2.02\times 10^{-13}$ & $ 1.08\times 10^{-11}$ \\
	16 & $7.77\times 10^{-15}$ & $1.06\times 10^{-13}$ & $1.60\times 10^{-13}$ & $ 7.90\times 10^{-12}$ 
	\end{tabular}
	\end{center}
	\caption{Maximum relative errors in a pipe flow simulation along lines parallel to the top wall. As the target points get closer to the wall, the integrands become more challenging to evaluate accurately, but special quadrature is activated giving high accuracy everywhere using a small amount of panels.  }\label{tab:pipe_flow}
\end{table}

\begin{figure}[!h]
	\begin{tabular}{c c}
		\includegraphics[width=0.45\textwidth]{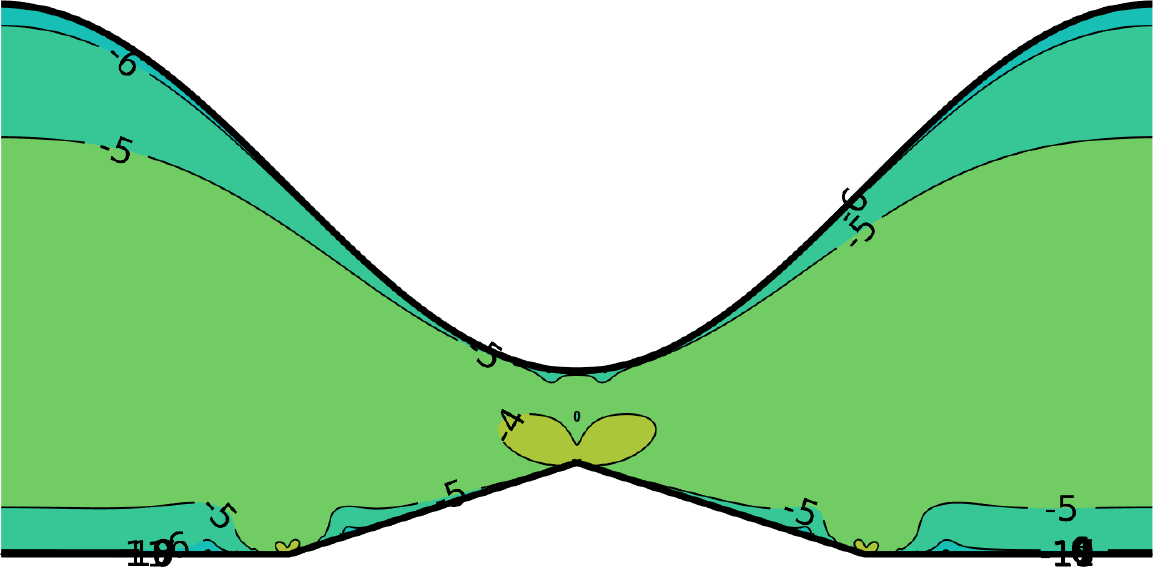} &
		\includegraphics[width=0.45\textwidth]{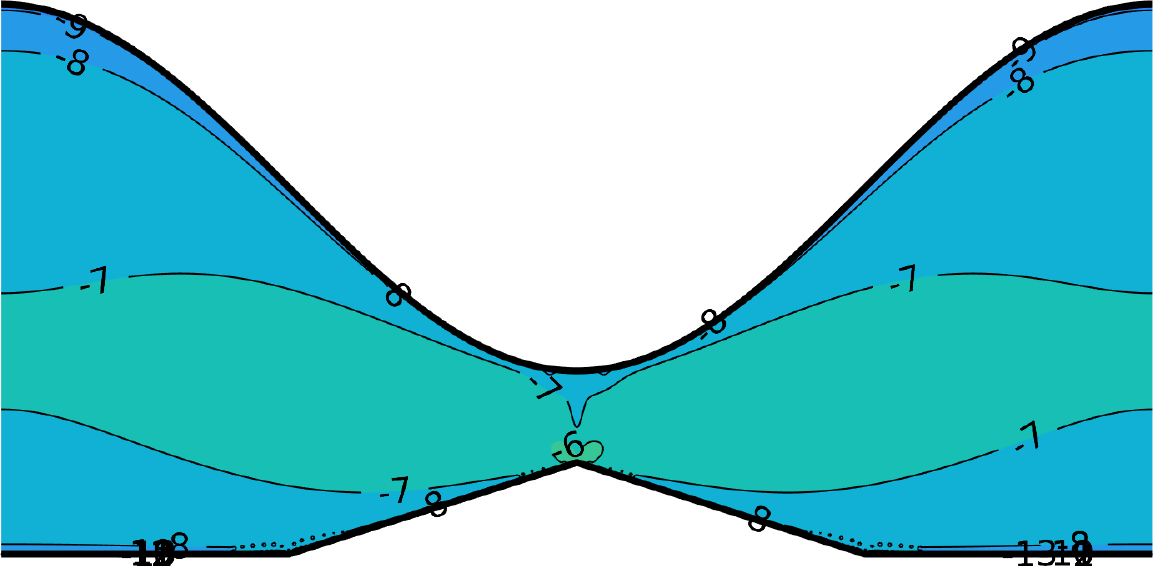} \\
		\includegraphics[width=0.45\textwidth]{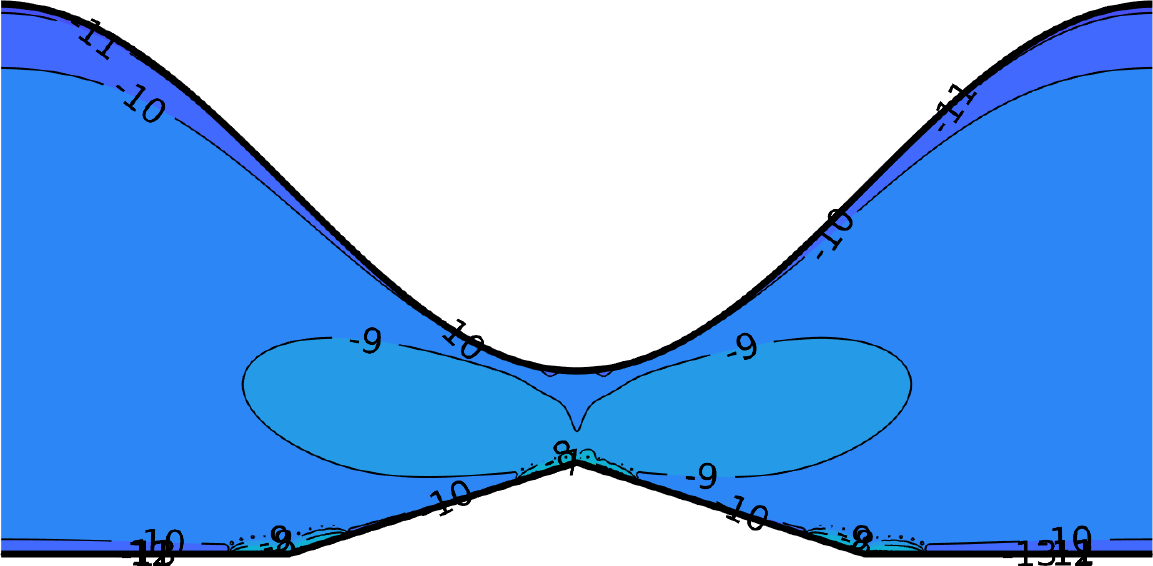} &
		\includegraphics[width=0.45\textwidth]{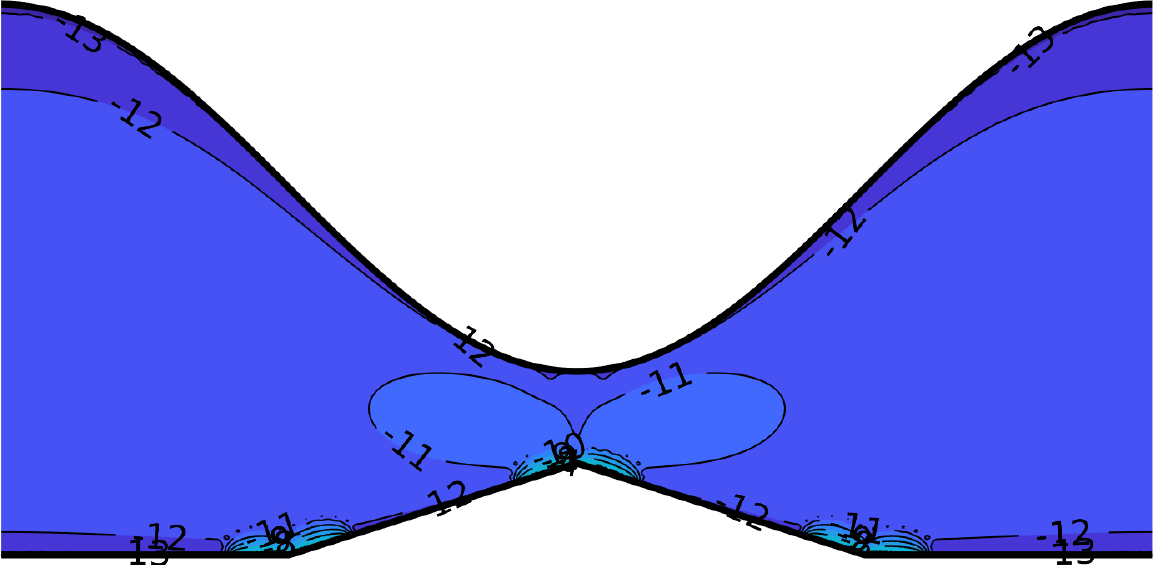}
	\end{tabular}
	\begin{center}
	\includegraphics[width=0.9\textwidth]{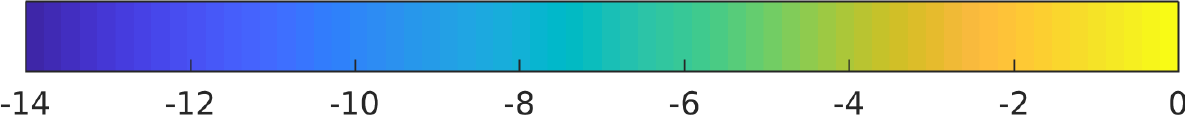}
	\end{center}
\caption{A shear flow $\uu = \langle x_2, 0\rangle$ is prescribed as a Dirichlet boundary condition on the top and bottom walls. If we take the pressure drop $\delta p$ to be 0, the BIE should recover the shear flow everywhere in the domain. Contour plots of the base 10 log of the relative error are shown for different values of $n_\sub$. Top: $n_\sub = 0,~10$, bottom: $n_\sub = 20,~30$. In every case each wall has been discretized using 20 panels with $n_q = 16$, meaning $N = 1280$ unknowns. Without any special treatment, the recovered solution differs quite a bit from the exact solution everywhere. With $n_\sub = 20$ the BIE can achieve 11 digits accuracy everywhere except very close to the corners. If additional accuracy is required near the corners, the fine density function can be recovered and used in conjunction with the special quadrature described earlier. }\label{fig:shear_corners}
\end{figure}

Another test is to prescribe Dirichlet boundary conditions on the solid walls. We will peform two types of tests: convergence towards a known, smooth solution, and self-convergence test towards an unknown solution. For the known smooth solution, we will prescribe a constant shear flow $\uu = \langle x_2, 0\rangle$ as boundary conditions on the top and bottom wall. The bottom wall is now only piecewise smooth. To recover a shear flow in the interior from the BIE solution, we must prescribe no pressure drop from inlet to outlet, i.e., $\langle \nabla p\rangle = 0$. We will use this simple example to test the RCIP algorithm for various values of $n_\sub$. Figure \ref{fig:shear_corners} shows the results. If we use the standard Nystr\"{o}m discretization without RCIP, we get quite large errors everywhere in the domain; the errors are not localized around the corner. For $n_\sub=30$ the RCIP algorithm allows us to achieve 11 digits accuracy everywhere in the domain, except for close to the corners. To gain additional accuracy near the corners, it is possible to recover the actual fine density and use it along with the special quadrature described in Section \ref{sec:post} \cite[Section 10]{Helsing2013b}.

\begin{figure}[!h]
\begin{tabular}{c c}
\includegraphics[width=0.45\textwidth]{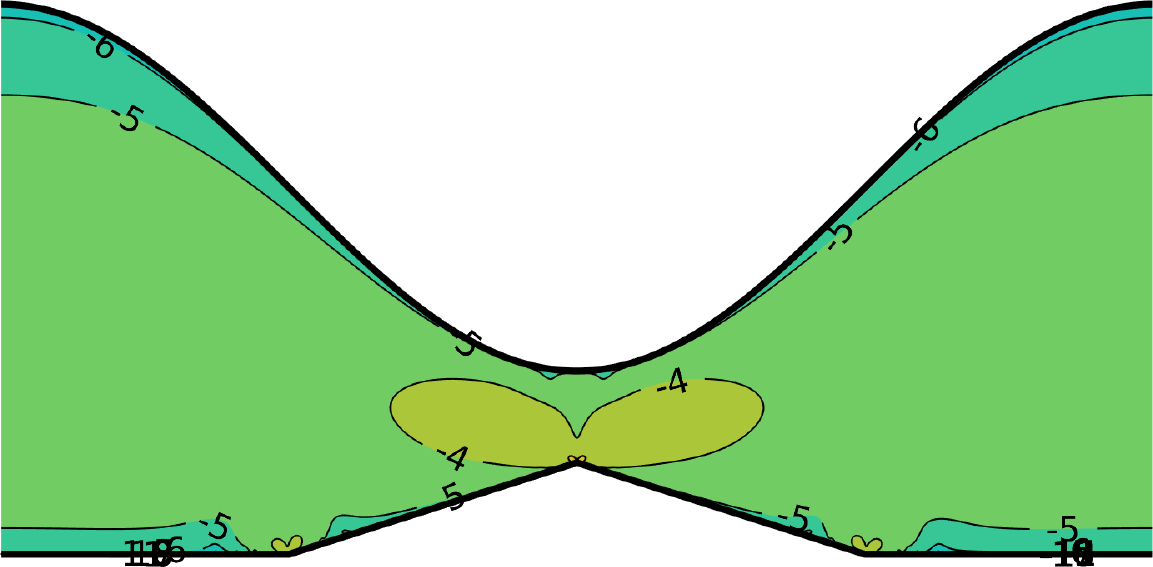} &
\includegraphics[width=0.45\textwidth]{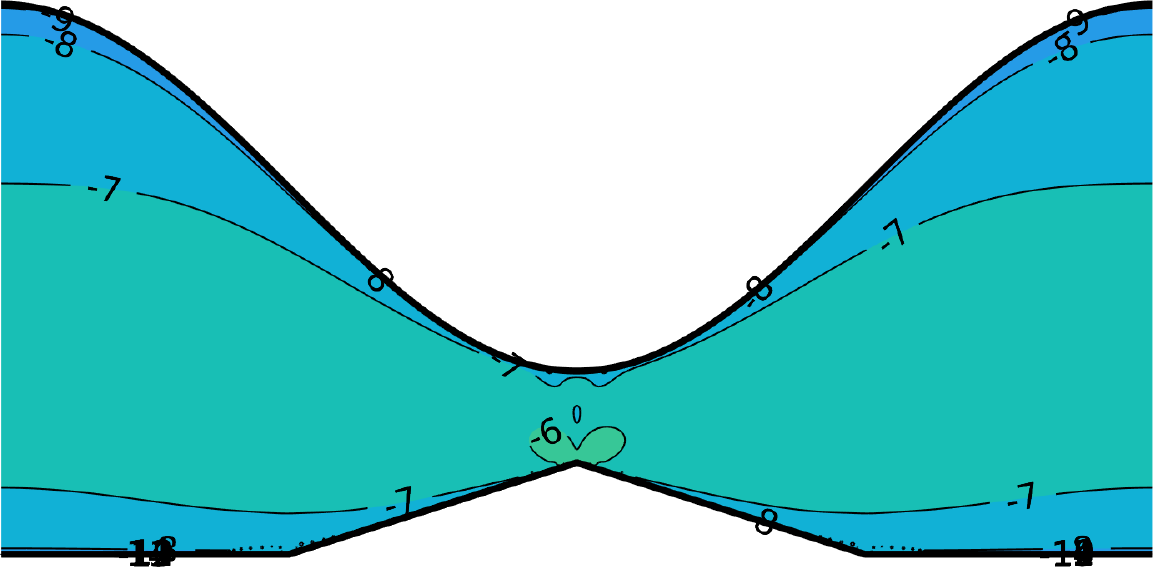} \\
\includegraphics[width=0.45\textwidth]{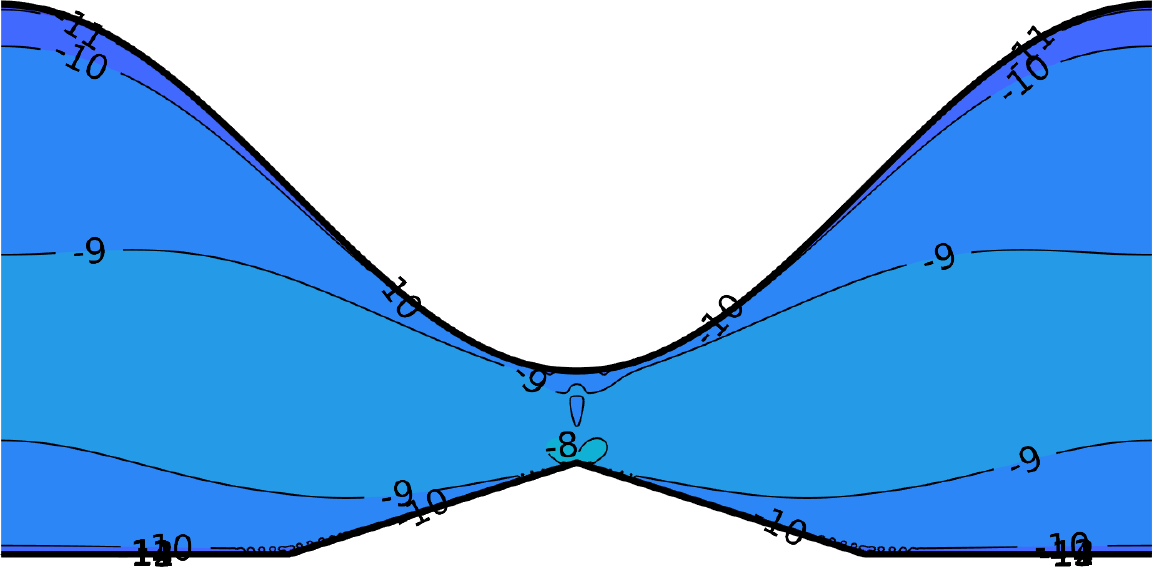} &
\includegraphics[width=0.45\textwidth]{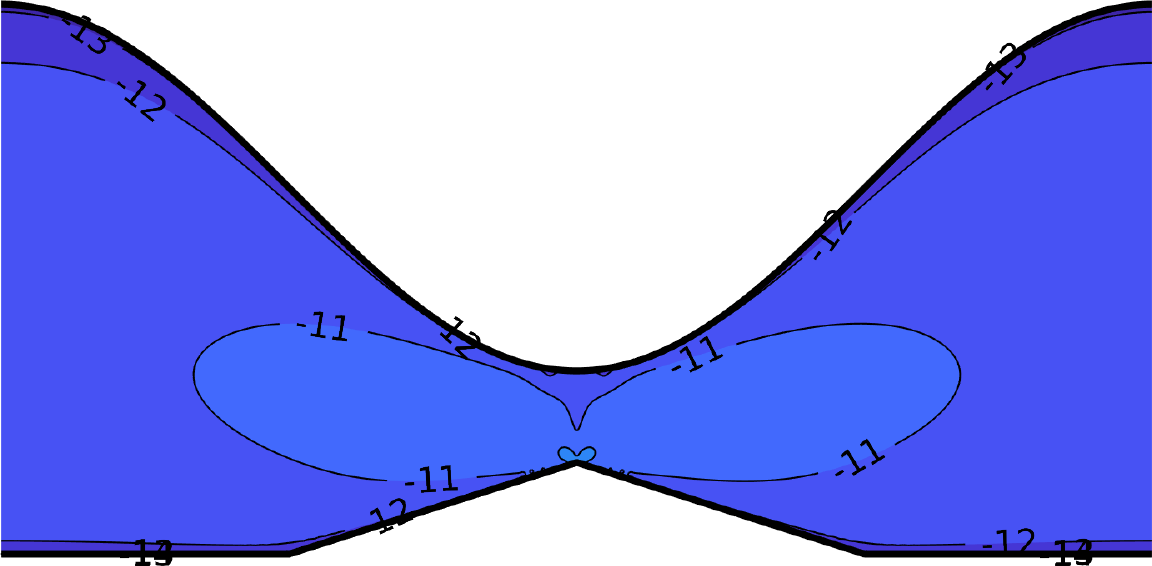}
\end{tabular}
	\begin{center}
	\includegraphics[width=0.9\textwidth]{colorbar.pdf}
	\end{center}
\caption{The same geometry as in Figure \ref{fig:shear_corners} is used to perform a self-convergence study. We prescribe an average pressure gradient $\langle p \rangle = 1$ and set $\uu = \mathbf{0}$ on the walls. We then compare the velocity field for different $n_\sub$ to one with $n_\sub = 50$. Top: $n_\sub = 0,~10$, bottom $n_\sub = 20,~30$.  Each wall has been discretized using 20 panels with $n_q = 16$, meaning $N = 1280$ unknowns for each simulation. As with Figure \ref{fig:shear_corners}, the contour plot shows the base 10 log of the relative error.}\label{fig:self_convergence}
\end{figure}

For the self-convergence study, we will use the same domain, but prescribe zero boundary conditions on the walls, and set $\langle \nabla p\rangle = 1$. We compare the computed velocity with $n_\sub = 0,~10,~ 20$, and 30, to one with $n_\sub = 50$.  As can be seen in Figure \ref{fig:self_convergence}, the velocity field using $n_\sub= 30$ matches the solution with $n_\sub = 50$ to 11 digits. 

\begin{figure}[!h]
\begin{tabular}{c c}
\includegraphics[width=0.45\textwidth]{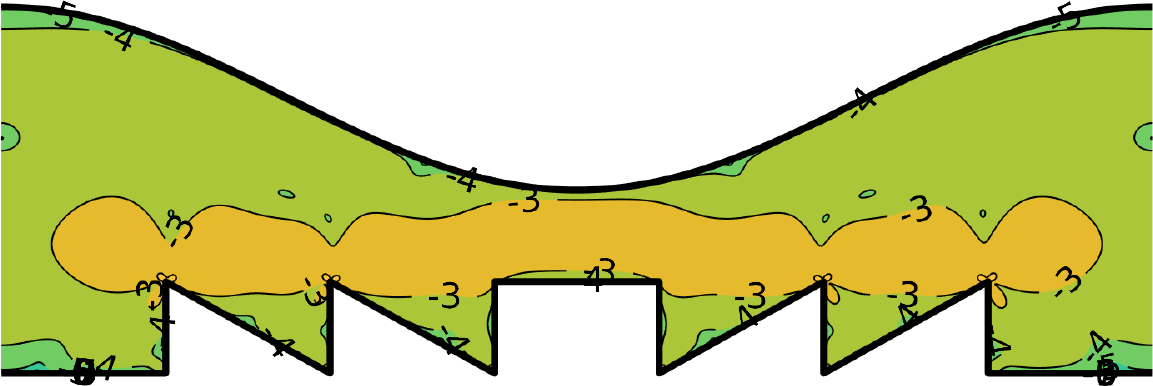} &
\includegraphics[width=0.45\textwidth]{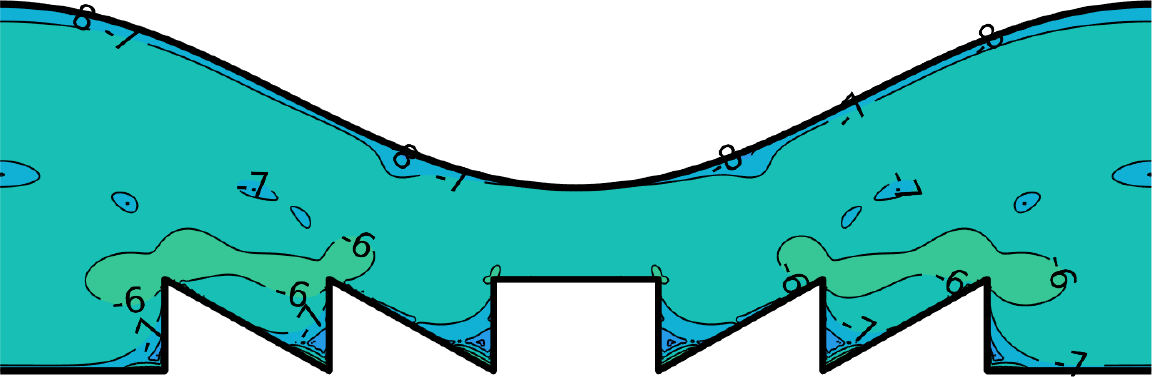} \\
\includegraphics[width=0.45\textwidth]{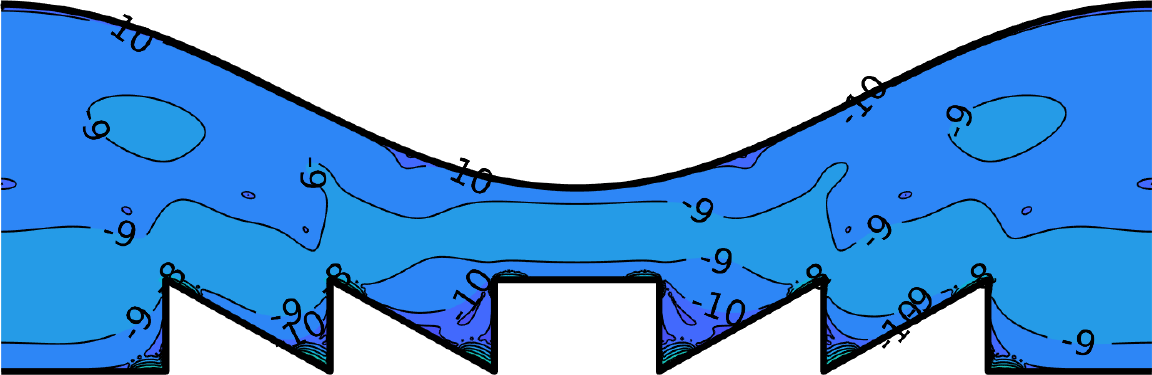}  &
\includegraphics[width=0.45\textwidth]{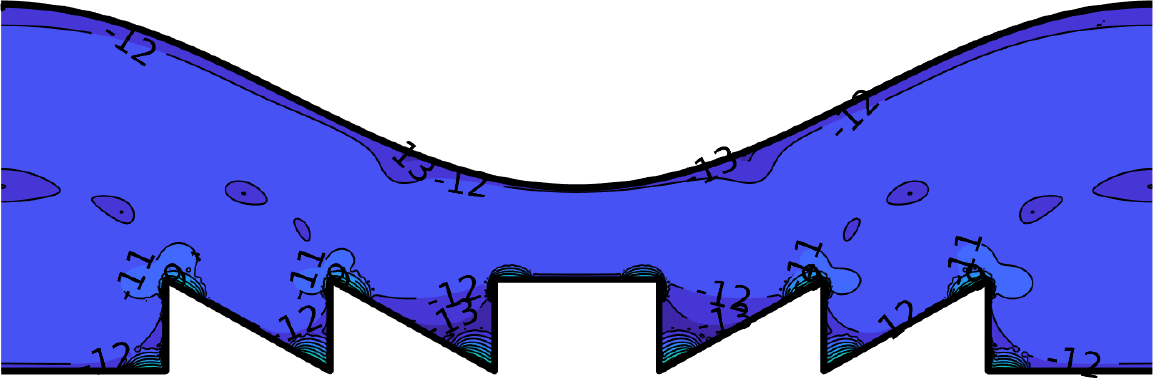}
\end{tabular}
	\begin{center}
	\includegraphics[width=0.9\textwidth]{colorbar.pdf}
	\end{center}
\caption{A shear flow $\uu = \langle x_2, 0\rangle$ is prescribed as a Dirichlet boundary condition on the top and bottom walls.  If we take the pressure drop $\delta p$ to be 0, the BIE should recover the shear flow everywhere in the domain. Note that the bottom wall contains several reentrant corners. Contour plots of the base 10 log of the relative error are shown for different values of $n_\sub$. Top: $n_\sub = 0,~20$, bottom: $n_\sub = 40,~60$. Both the top and bottom walls are discretized with 65 panels with $n_q=16$, meaning 4160 unknowns. For comparison, a locally refined solution with $n_\sub = 60$ for this problem would require over 50,000 unknowns. }\label{fig:shear_corners2}
\end{figure}

The number of refinements needed to achieve a desired accuracy depends heavily on the geometry. Reentrant corners in particular require more refinement. For example Figure \ref{fig:shear_corners2} shows an example where $n_\sub = 20$ achieves only a maximum of 6 digits of accuracy. Using $n_\sub = 60$ we are able to acheive 11 digits accuracy in the interior. At the finest $n_\sub$ we are placing quadrature points on panels of size $10^{-19}$. Using a local coordinate system centered at zero is essential in this case.

\subsection{Remarks on timings} 

Each subdivision in the RCIP algorithm involves inverting an $\mathbf{R}$ matrix in \eqref{eq:recursion}. This matrix will be of the size $8n_q\times 8n_q$; all our simulations use  $n_q = 16$, meaning $\mathbf{R}$ will be $128\times 128$. This is not a very cheap procedure, and for small problems this cost can indeed dominate the cost of solving the linear system to find $\tilde{\qq}$. Using Matlab 2017a on a desktop with 32 GB of RAM and an Intel i7 processor each subdivision step takes approximately 1.6 seconds (this is  indpendent of the problem under consideration). For the problem in Figures \ref{fig:shear_corners} and \ref{fig:self_convergence}  with $n_\sub=10$, assembling the RCIP matrices (3 of them) takes around 48 seconds. By comparison  solving the linear system of size $1280\times 1280$ using GMRES takes around 0.7 seconds, and postprocessing the solution at 25,000 points (and applying special quadrature where needed) takes around 4 seconds. 

The cost of a single matrix-vector product is slightly higher as well using RCIP. Without RCIP, for each GMRES iteration it is necessary to compute the product $\mathbf{K}\qq$, ususally using some kind of fast summation method. As discussed in Section \ref{sec:multiplication}, when using RCIP, each matrix-vector product requires the additional computation of $\hat{\RR}\tilde{\qq}$, as well as the product $\mathbf{K}^*(\hat{\RR}\tilde{\qq})$. As mentioned, the first of those products can be done directly, since $\hat{\RR}$ is a block diagonal matrix, with one $8n_q\times 8n_q$ block per corner. The second of these products is actually just another $8n_q\times 8 n_q$ product per corner. This product is effectively zeroing out the entries in the matrix $\mathbf{K}$ corresponding to entries where both the source and target points are close to corners. Therefore it should be done using the same fast method as the full matrix-vector product $\mathbf{K}(\hat{\RR}\tilde{\qq})$.  In the example corresponding to Figures \ref{fig:shear_corners} and \ref{fig:self_convergence}, one full matrix-vector product takes 0.035 seconds without RCIP and 0.038 seconds with RCIP, an increase of about 8\%. This percentage increase will depend on the $n_q$, the number of corners, as well as the relative size of $\Gamma^*$ to $\Gamma$, but \emph{not} on $n_\sub$.

It should be noted, that without RCIP the cost grows as $\mathcal{O}(N\log N)$ with local refinement. Here $N$ is the number of unknowns, and grows by $4 n_q$ per corner for each subdivision. In the examples shown in Figures \ref{fig:shear_corners} and \ref{fig:self_convergence} this would correspond to solving a linear system of size $1280 + 4(16)(10)(3) = 3200$; for $n_\sub = 50$, the linear system would be of the size 10,880. In addition, the condition number of the refined matrix would grow with the number of quadrature points, and require more GMRES iterations to solve to a desired accuracy. For time dependent problems as in Section \ref{sec:drops}, the larger system would have to be solved at every time step. In this case RCIP has a clear advantage, since as long as the wall geometry remains fixed, the $\RR$ matrices need only be assembled once and a smaller linear system can be solved at each time step.

\section{Adding Viscous Drops}\label{sec:drops}

Earlier work \cite{Palsson2019b} has looked at modelling the movement of drops inside confined periodic geometries.  We now demonstrate the robustness and usefulness of RCIP by extending the method in \cite{Palsson2019b} to model the movement of drops near sharp interfaces. A drop is a packet of fluid that is does not mix with the fluid surrounding it. Surface tension forces prevent the mixing of the drop with the bulk fluid. Both the fluid inside each drop, and the bulk fluid, satisfy the incompressible Stokes equations,

\begin{alignat*}{2}
	-\mu_0 \Delta \uu(\xx) + \nabla p(\xx) &= \mathbf{0},\qquad  &&\xx\in\Omega_0,\\
	-\mu_\ell \Delta\uu(\xx) + \nabla p(\xx) &= \mathbf{0},\qquad &&\xx\in\Omega_\ell,~\ell=1,\hdots, N_d,\\
	\nabla\cdot\uu &= 0,\qquad &&\xx\in\bigcup\limits_{\ell=0}^{N_d}\Omega_\ell,
\end{alignat*}
where $\mu_\ell$ denotes the viscosity inside the region bounded by $\Gamma_\ell$. For convenience we will define the viscosity ratio $\lambda_\ell = \mu_\ell / \mu_0$. As $\lambda$ increases the drop behaves more like a rigid particle.

\begin{figure}[!h]
\begin{center}
\includegraphics[width=0.9\textwidth]{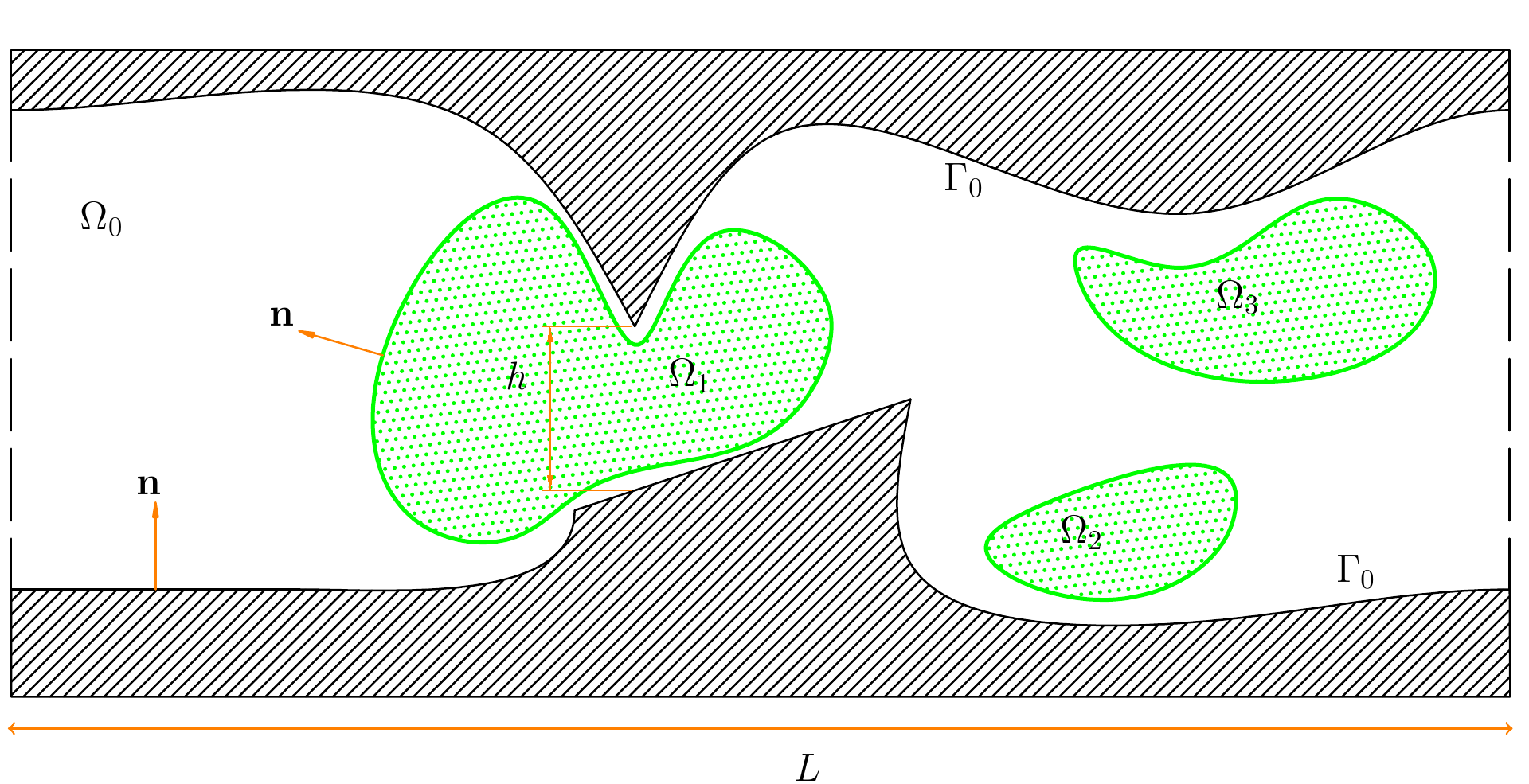}
\end{center}
\caption{Sketch of a problem involving deformable drops. The region $\Omega_0$ is the bulk fluid and $\Gamma_0$ are solid walls. The regions $\Omega_1$, $\Omega_2$, and $\Omega_3$ are drops, with interfaces $\Gamma_1$, $\Gamma_2$ and $\Gamma_3$ respectively. On each boundary the normal vector $\nn$ points into the bulk fluid. The problem is periodic in the horizontal direction.}\label{fig:drop_sketch}
\end{figure}

At the interfaces the velocity is continuous; however, in general the surface forces acting on the interface from the inside and the outside of the drop are not equal. We will denote the jump in the normal force across interface $\ell$ as $\delta \ff_\ell(\xx)$. This jump is related to the curvature $\kappa_\ell(\xx)$, and the surface tension $\sigma_\ell(\xx)$ by \cite{Pozrikidis1992} [Chapter 5],
\[\delta \mathbf{f}_\ell(\xx) =: (-p_0 -p_\ell)\nn_\ell + 2(\mu_0 \mathbf{e}_0 - \mu_k\mathbf{e}_\ell)\cdot \nn_\ell  = \sigma_\ell(\xx) \kappa_\ell(\xx) \nn(\xx) - \nabla_s \sigma_\ell(\xx),\]
where $\mathbf{e}_\ell$ denotes the rate of strain tensor in $\Omega_\ell$ and $\nabla_s$ is the gradient along the interface. 

To nondimensionalize this problem, we will define $h$ to be the minimum height of the channel. We would like the maximum velocity of an empty channel to be 1. To do this, we will introduce a pressure scale $h\langle \nabla p\rangle/8$,  a length scale $h$, a velocity scale $h^2\langle \nabla p\rangle /(8 \mu_0)$, and a surface tension scale $\sigma_0$ to obtain the full nondimensional form of the problem,
\begin{subequations}
	\begin{alignat*}{2}
	-\Delta \uu(\xx) + \nabla p(\xx) &= \mathbf{0},\qquad  &&\xx\in\Omega_0,\\
	-\lambda_\ell \Delta\uu(\xx) + \nabla p(\xx) &= \mathbf{0},\qquad &&\xx\in\Omega_\ell,~\ell=1,\hdots, N_d,\\
	\nabla\cdot\uu &= 0,\qquad &&\xx\in\bigcup\limits_{\ell=0}^{N_d}\Omega_\ell,\\
	\delta\ff_\ell &= \frac{1}{\text{Lp}}(\sigma_\ell \kappa_\ell \nn_\ell - \nabla_s\sigma_\ell), \qquad && \xx\in\Gamma_\ell,\\
	\uu(\xx) &=\uu(\xx + L),\\
	\langle \nabla p\rangle &= (8, 0),
\end{alignat*}
\end{subequations}
where the Laplace number $\text{Lp}$ is the dimensionless quantity given by $h^2\langle\nabla p\rangle /(8\sigma_0)$ \cite{Kunes2012}.  The time scale is now $h / U =8 \mu_0/(h \langle \nabla p\rangle)$. The Laplace number serves the same purpose as the Capillary number in other drop simulations \cite{Palsson2019b}, but for pressure driven flows it is more natural to nondimensionalize according to the pressure gradient as opposed to a maximum velocity. In \cite{Palsson2019a,Palsson2019b} chemical surface active agents change the surface tension of the drops; however, for the remainder of this paper we will assume a constant surface tension, i.~e. $\nabla_s \sigma = 0$. Allowing for non constant surface tension would not impact the RCIP algorithm in any way. 

A BIE formulation \cite{Zinchenko2006} for the problem is given by
\begin{equation}\label{eq:zinchenko} \uu(\xx) = \sum\limits_{\ell = 1}^{N_d} \left(\mathbb{S}^P_{\Gamma_\ell} [\delta\ff](\xx) + (\lambda_\ell - 1)\mathbb{D}^P_{\Gamma_\ell}[\uu](\xx)\right) + \mathbb{K}^P[\qq](\xx) + \langle \uu \rangle,\end{equation}
    where $\mathbb{S}^P_{\Gamma_\ell}$ and $\mathbb{D}^P_{\Gamma_\ell}$  are the periodic single- and double-layer operators defined on $\Gamma_\ell$.
    
Taking the limit as $\xx\to\xx_0\in \Gamma_0$,
\[ -\qq(\xx_0) + \mathbb{K}^P[\qq](\xx_0) + \sum\limits_{\ell = 1}^{N_d}  (\lambda_\ell - 1)\mathbb{D}^P_{\Gamma_\ell}[\uu](\xx_0) + \langle \uu \rangle = \mathbf{g}(\xx_0) - \sum\limits_{\ell=1}^{N_d}\mathbb{S}^P_{\Gamma_\ell} [\delta\ff](\xx_0), \quad \xx_0\in\Gamma_0,\]
and taking the limit $\xx\to\xx_0 \in\Gamma_m$,
\begin{align*}
\uu(\xx_0) -2\sum\limits_{\ell = 1}^{N_d} \frac{\lambda_\ell - 1}{\lambda_m + 1} \mathbb{D}^P[\uu](\xx_0) &- \frac{2}{\lambda_m + 1}\left(\mathbb{K}^P[\qq](\xx_0) + \langle\uu\rangle\right) \\& = \sum\limits_{\ell=1}^{N_d} \frac{2}{\lambda_m + 1}\mathbb{S}^P_{\Gamma_\ell}[\delta\ff](\xx_0) \qquad \xx_0\in\Gamma_m.
\end{align*}

As in Section \ref{sec:periodicity} we close the system by relating the density function $\qq$, still defined only on the solid walls $\Gamma_0$, to the (nomdimensionalized) average pressure gradient,
\[ \eta\int_{\Gamma_0} \qq ~\text{d}\Gamma_0 = \frac{1}{L^2}\begin{pmatrix} 8 \\ 0\end{pmatrix}.\]

By inspecting \eqref{eq:zinchenko}, we see that if $\lambda_\ell = 1,~\ell = 1,\cdots,N_d$, then the velocity on the drop interfaces only appears on one side of the equation. This means that after solving for $\qq$ on the solid walls, $\uu$ can be computed as a post-processing step. In general; hower, both $\uu$ and $\qq$ need to be solved for simultaneously. The drop boundaries are discretized in exactly the same way as the solid walls, and the same close evaluation scheme is used to compute both $\mathbb{D}_{\Gamma_\ell}^P[\uu]$ and $\mathbb{S}_{\Gamma_\ell}^p[\delta\ff]$ for drops that are near each other, or near a solid wall.  After computing $\uu(\xx)$ on the drop interfaces, the drops move and deform according to the velocity,
\[ \frac{\text{d}\xx}{\text{d}t} = \uu(\xx), \qquad \xx \in \bigcup\limits_{\ell = 1}^{N_d}\Gamma_\ell.\]
To evaluate this system of ODEs we will use the fourth order adaptive time stepper described in \cite{Kennedy2003}. Further details on time stepping, including a way to maintain consistent acrclength spacing as the drop perimeter changes, can be found in \cite{Palsson2019a,Palsson2019b}.

It is important to note that since the solid walls are not moving, the RCIP matrices $\RR_{j}$, $j=1,\hdots, n_c$ need only be computed at the start of the simulation. This means that after this precomputation the actual time needed each time step to solve the linear system is independent of $n_\sub$. 

\subsection{Numerical results involving visous drops}

\begin{figure}[!h]
\begin{tabular}{c c c}
$t = 0$ & $t = 0.5$ & $t=1$\\
\includegraphics[width=0.3\textwidth]{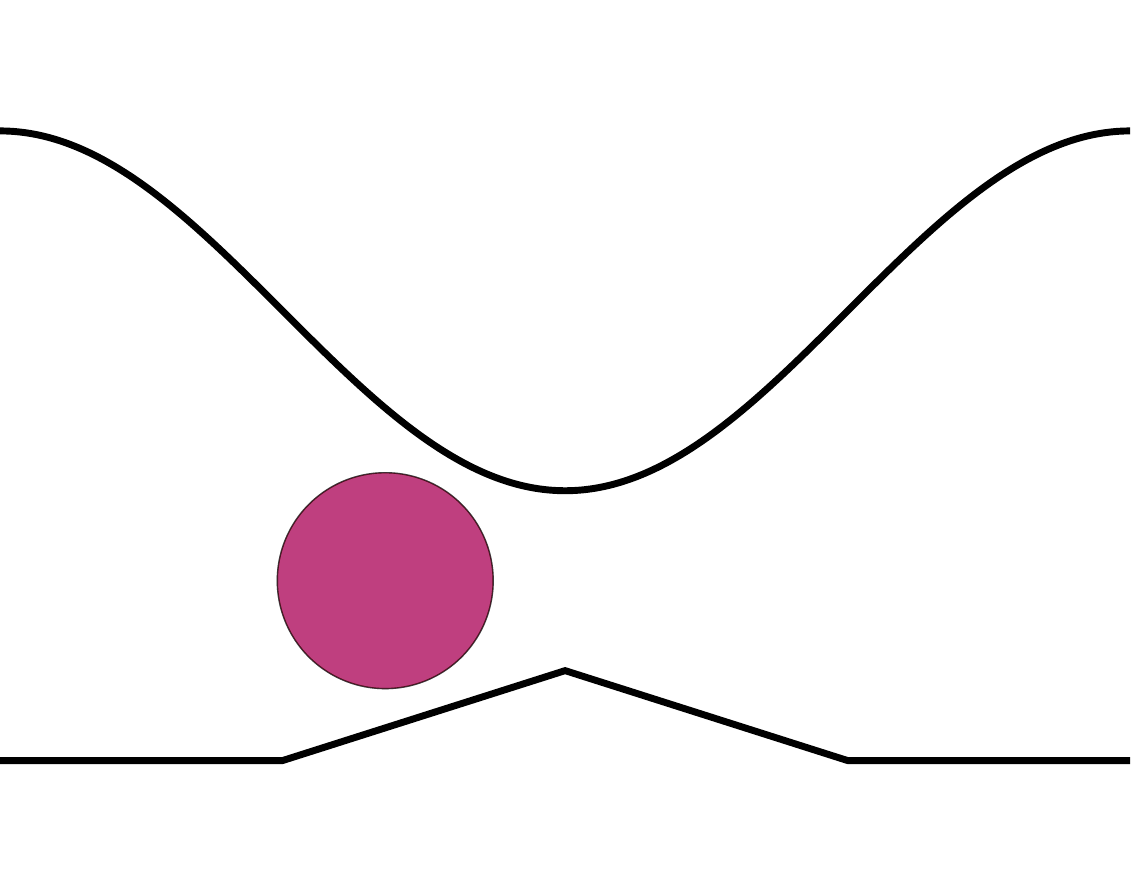} &
\includegraphics[width=0.3\textwidth]{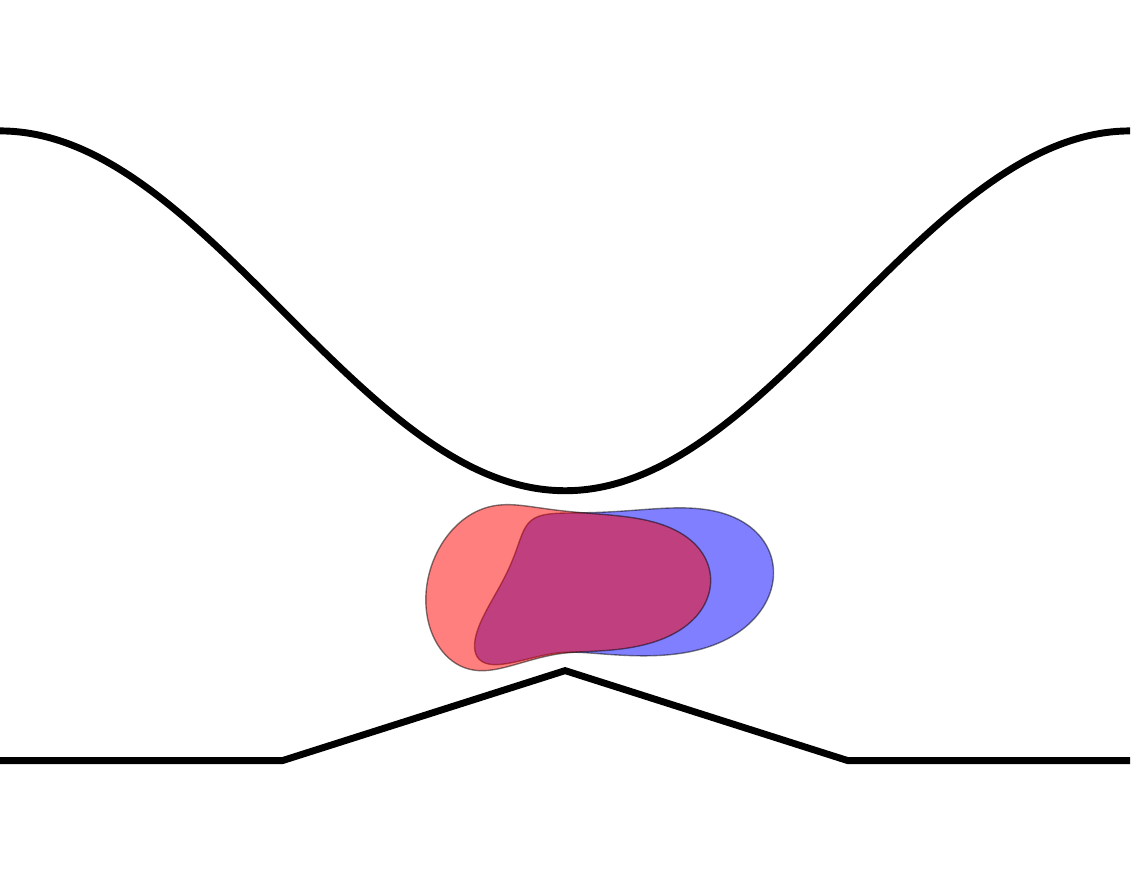} &
\includegraphics[width=0.3\textwidth]{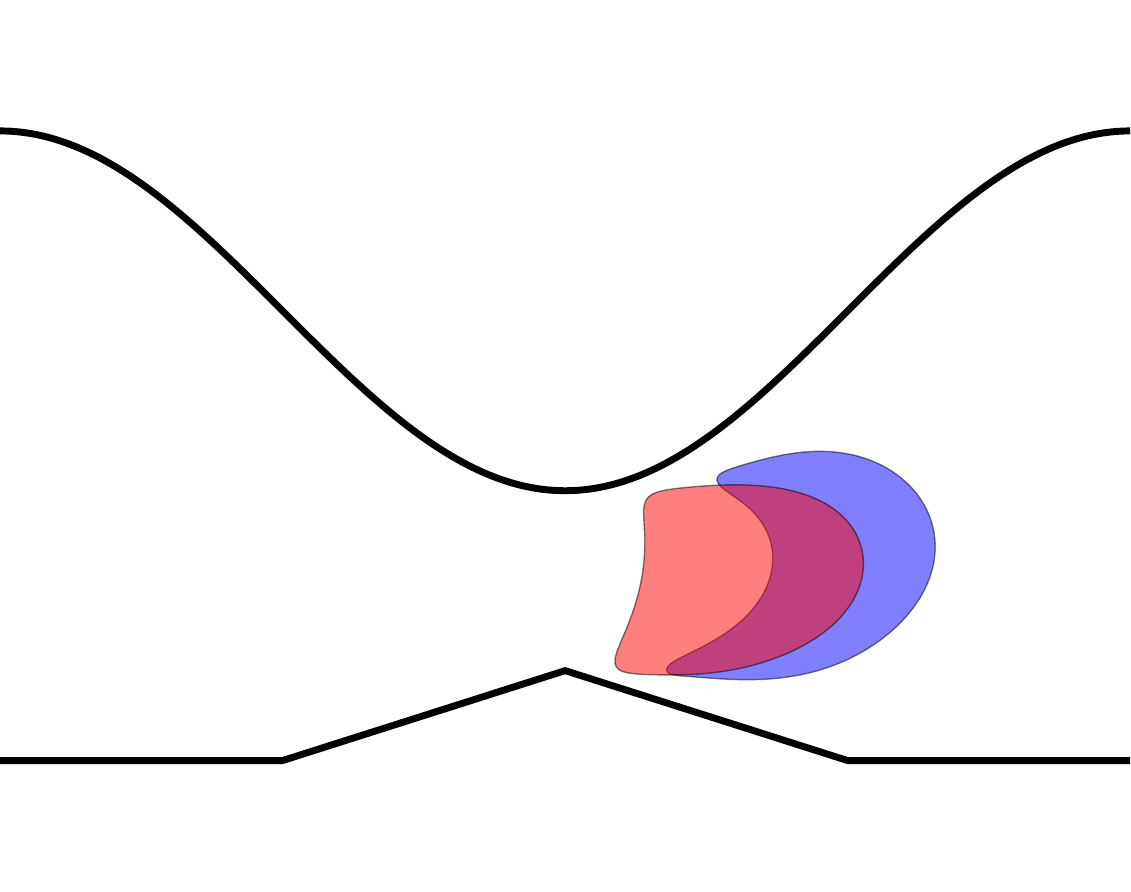}
\end{tabular}
\caption{Snapshots of a drop moving inside a channel containing a corner. A pressure drop is applied from left to right. By varying the number of points on the drop and the wall we can perform a numerical convergence study, the results of which are shown in Figure \ref{fig:convergence_plots}. We have run two simulations, one with $\lambda=1$ (blue drop), and another with $\lambda=5$ (red drop).  As expected, the drop with $\lambda=5$ deforms less than the $\lambda=1$ drop. $Lp=1$ for all simulations. Note that for presentation purposes the red and blue drops are overlayed on top of each other and appear to be purple where they overlap. }\label{fig:convergence_snapshots}
\end{figure}

To illustrate the performace of our method, we will consider two different examples. First, figure \ref{fig:convergence_snapshots} shows snapshots of a simulation of a single drop passing through a narrow piecewise continuous constriction. We begin by creating a high resolution reference solution with RCIP for $\lambda=1$ and $\lambda=5$. The reference solution will be discretized using $n_{\text{pan}}= 140$ on both the solid walls and the drop, and $n_\sub=20$. The time stepping tolerance for the high resolution simulation is $10^{-10}$. 

To compute an error, we will look at the $\ell_\infty$ difference of a numerical solution with the computed reference solution. As described in \cite{Palsson2019b}, when updating the positions of the drops it is advantageous to work on a uniform grid with parameter spacing $\Delta s$, instead of the panel Gauss-Legendre points. Using a uniform grid also allows us to easily upsample the non-reference solutions using an FFT to obtain a solution at the same discretization points as the reference solution. 

\begin{figure}[!h]
\begin{center}
\includegraphics[width=0.8\textwidth]{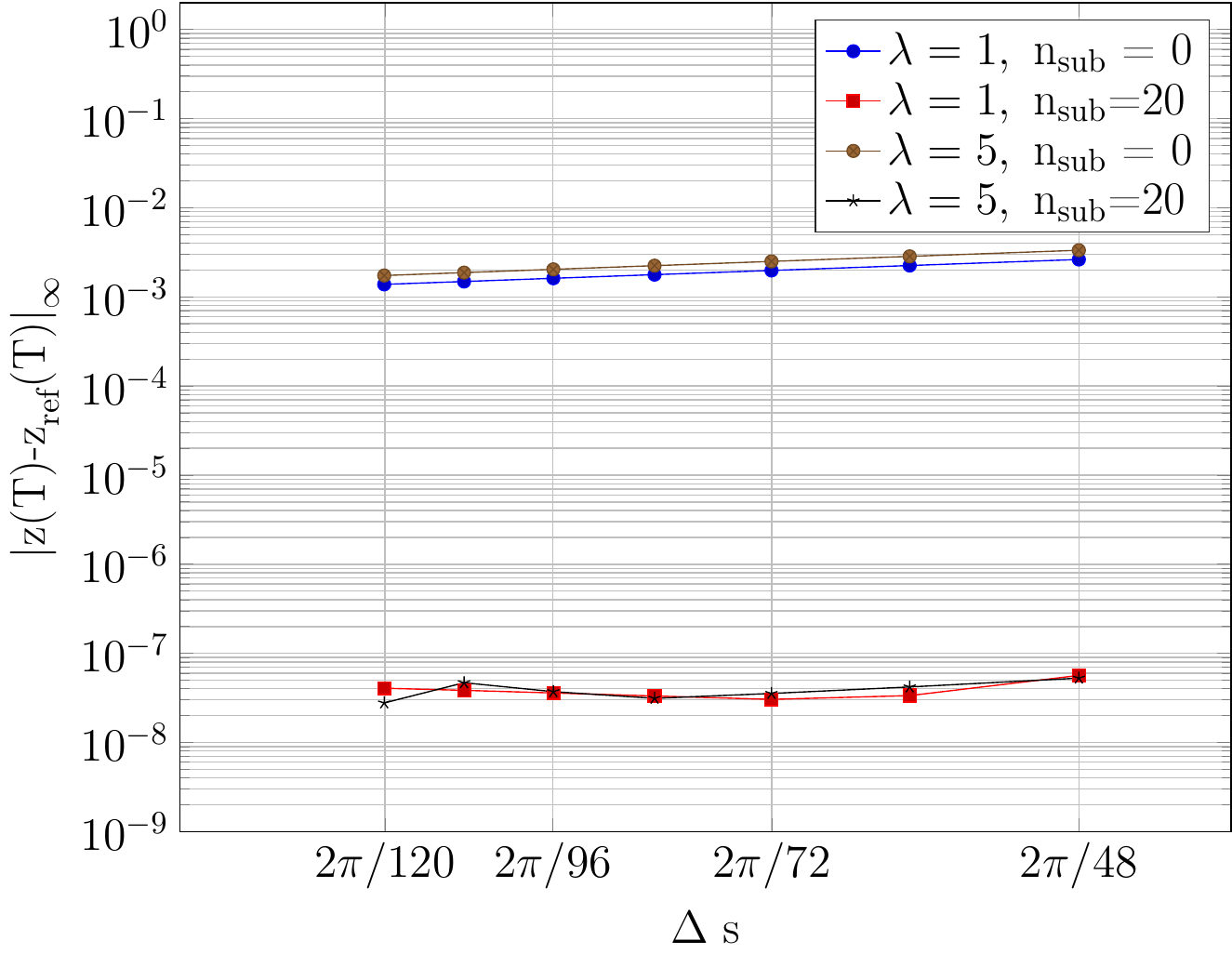} 
\end{center}
\caption{Spatial error study with and without RCIP for the simulation shown in Figure \ref{fig:convergence_snapshots}. Using RCIP gives much lower errors, and these are near the time stepping tolerance of $10^{-8}$. $Lp=1$ for all simulations. }\label{fig:convergence_plots}
\end{figure}

As can be seen in Figure \ref{fig:convergence_plots}, without RCIP, as we refine the number of points on both the drop and the wall, the numerical solutions for both $\lambda=1$ and $\lambda=5$ converge very slowly towards the high resolution reference solution. When we use RCIP with $n_\sub=20$, the numerical solutions are much more accurate; around the level of the time stepping tolerance for the low resolution simulations of $10^{-8}$.

\begin{figure}[!h]
\begin{center}
\includegraphics[width = 0.8\textwidth]{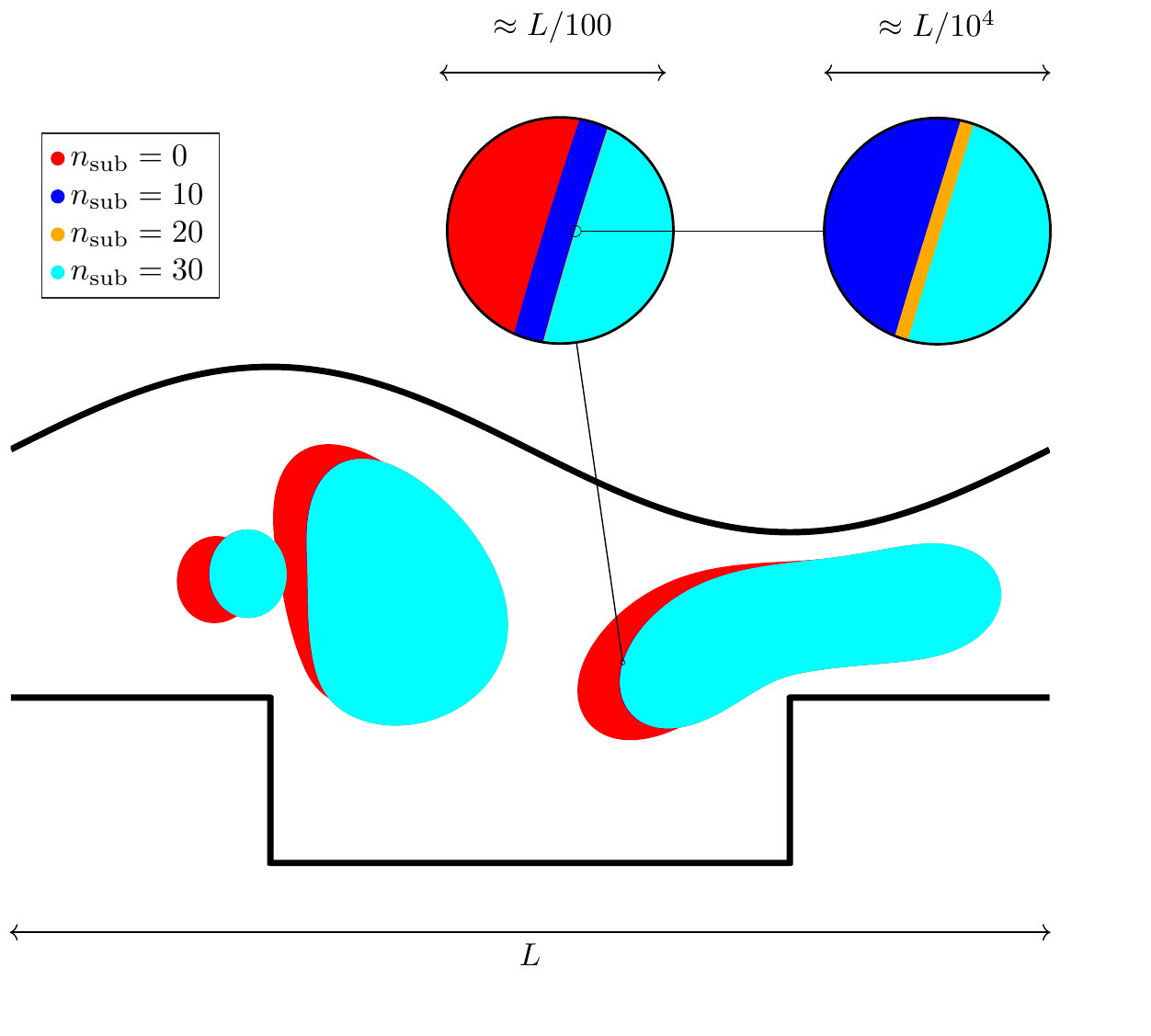}
\end{center}
\caption{Snapshot of three drops of different sizes in a periodic channel. Four different simulations are run starting from identical initial conditions, each with a different $n_\sub$. For the red drops, no RCIP is used;  $n_\sub = 10$ for blue drops; $n_\sub = 20$ for orange drops; $n_\sub=30$ for cyan drops. The area error as a function of time is shown in Figure \ref{fig:squeezing_error}.  There is a visible difference between the $n_\sub = 0$ and all the other simulations, thus demonstrating the necessity of the special corner treatment. As $n_\sub$ increases the magnified regions demonstrate that the drops converge to the same solution (the cyan drops). For all cases  $\lambda=1$ and $Lp=0.2$ for all the drops. The top wall is discretized with 20 panels, the bottom with 30 panels, and each drop with 60 panels. The time stepping tolerance is set to $10^{-10}$. }\label{fig:squeezing_snapshots}
\end{figure}

As a second example, we will investigate the movement of multiple drops of different sizes confined in a periodic channel containing several corners. A snapshot of such a simulation at $t=6$ (after all the drops have passed at least one periodic length) is shown in Figure \ref{fig:squeezing_snapshots}. As the snapshot makes clear, there is a clear, visible difference in the modelled drops depending on whether or not we use RCIP. In particular, the difference between no RCIP and $n_\sub = 10$ is quite large, whereas the difference between the $n_\sub=20$ and $n_\sub=30$ is hardly noticable. 

\begin{figure}[!h]
\begin{center}
	\includegraphics[width = 0.7\textwidth]{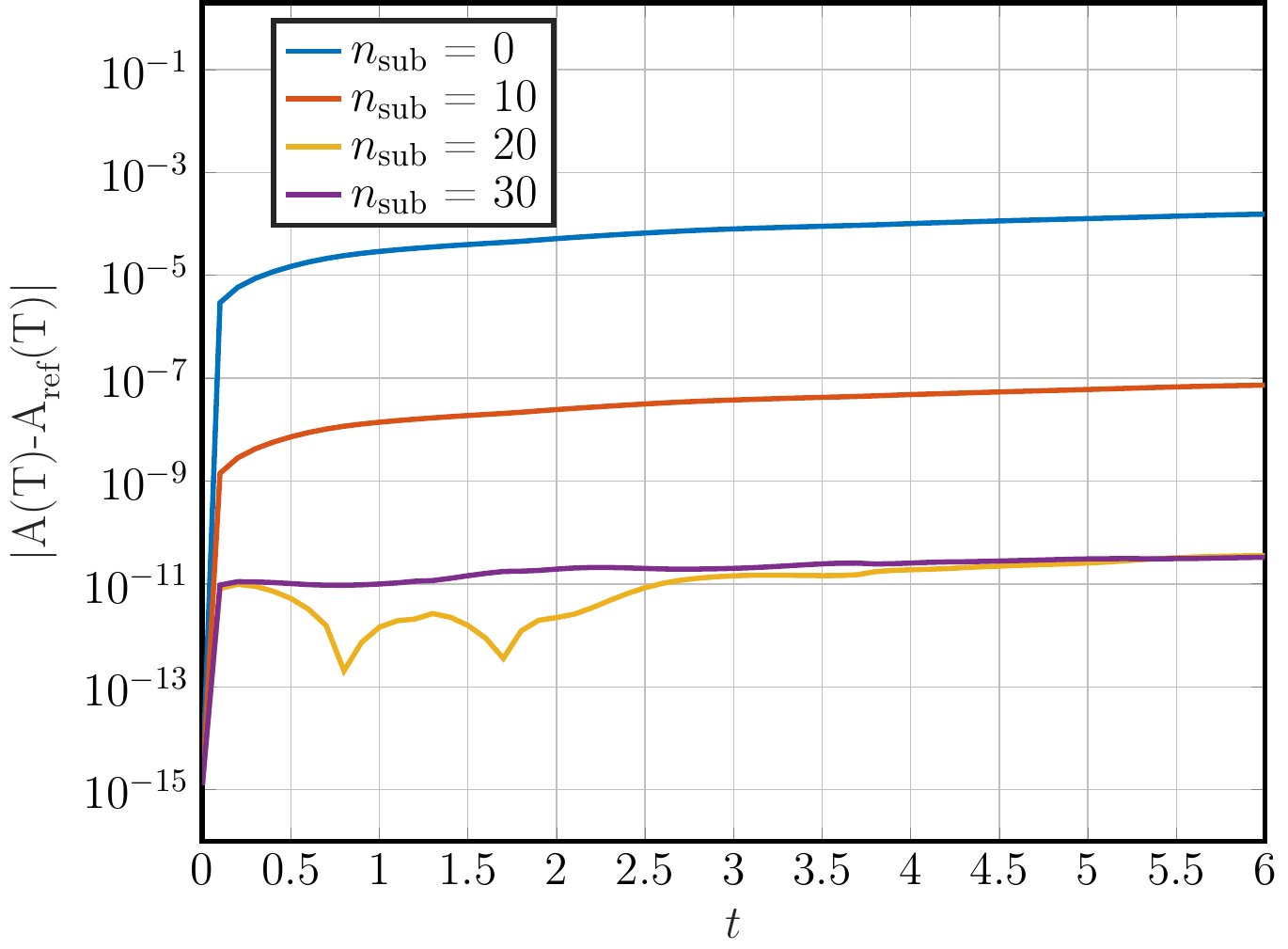}
\end{center}
\caption{Area conservation errors for the simulation shown in Figure \ref{fig:squeezing_snapshots} as a function of time. The error eventually plateaus around $10^{-4}$ without RCIP, but remains below $10^{-10}$ when RCIP with $n_\sub=20$ is used. The time stepping tolerance is set to $10^{-10}$. }\label{fig:squeezing_error}
\end{figure}

Further proof of the necessity of the proper handling of corners is shown in Figure \ref{fig:squeezing_error}. Since the fluid inside the drops is incompressible, the area of each drop should be conserved. Note that this is not enforced explicitly, nor is area conservation used in the criteria for the temporal adaptivity. Without RCIP, the area after each drop has passed one periodic channel length is conserved only up to around $10^{-4}$. When applying RCIP with $n_\sub =10$, it is below $10^{-7}$, and when applying $n_\sub=20$, the area error is at or below the time stepping tolerance. 

\section{Conclusion}

Domains with sharp corners pose challenges for boundary integral methods. The layer density defined on the boundary becomes weakly singular at the corners, and therefore standard quadrature rules cannot be accurately applied. This destroys the accuracy of the post-processed solution everywhere in the domain. We have demonstrated how a technique known as recursively compressed inverse preconditioning (RCIP) can be used to accurately solve the Stokes equations on two-dimensional domains with corners. This method requires a small amount of precomputation, however, it does not add any additional unknowns to the problem, nor does it increase the condition number of the linear system. Numerical experiments have demonstrated its robustness and usefulness for postprocessing the velocity anywhere in the domain. Additionally we have shown that it can be added to an existing drop model \cite{Palsson2019b} to accurately model the movement of drops near corners.

Future work could include extending this model to three dimensions. A boundary integral method for three-dimensional drops in free space has been developed  \cite{Sorgentone2018,Sorgentone2019}. In three dimensions Lipschitz domains admit both corners and edges, so the RCIP algorithm described in Section \ref{sec:rcip} must be further developed. In \cite{Helsing2013a} such an approach is used to compute the capacitance of a cube. 

 \begin{acknowledgements}

The authors would like to acknowledge support   the Knut and Alice Wallenberg Foundation (award number KAW 2017.0410) as well as by the Göran Gustafsson Foundation for Research in Nature and Medicine. 

\end{acknowledgements}
\bibliographystyle{plain}
\bibliography{bibliography}

\end{document}